\documentclass[a4paper]{article}
\usepackage[margin=2cm]{geometry}
\usepackage{graphicx}
\usepackage{xcolor}
\usepackage{amsmath}
\usepackage{mathtools}
\usepackage{bm}
\usepackage{xcolor}
\usepackage{amsfonts}
\usepackage{booktabs}
\usepackage{array}
\usepackage{amsthm}
\newtheorem{definition}{Definition}
\newtheorem{theorem}{Theorem}
\newtheorem{proposition}{Proposition}
\newtheorem{corollary}{Corollary}
\usepackage{hyperref}
\makeatletter
\newcommand{\doubletilde}[1]{{%
  \mathpalette\double@tilde{#1}%
}}
\newcommand{\double@tilde}[2]{%
  \sbox\z@{$\m@th#1\tilde{#2}$}%
  \ht\z@=.9\ht\z@
  \tilde{\box\z@}%
}
\makeatother
\usepackage{mathdots}

\title{Preconditioners for Fractional Diffusion Equations Based on the Spectral Symbol}
\author{Nikos Barakitis \\ Department of Informatics, Athens University of Economics and Business,\\ 76 Patission Str., GR-10434 Athens, Greece \and Sven-Erik Ekstr\"om \\ Department of Information Technology, Division of Scientific Computing,\\ Uppsala University, L\"agerhyddsv\"agen 1, 751 05 Uppsala, Sweden\and Paris Vassalos\thanks{pvassal@aueb.gr} \\ Department of Informatics, Athens University of Economics and Business,\\ 76 Patission Str., GR-10434 Athens, Greece}
\begin{document}

\maketitle

\abstract{
It is well known that the discretization of fractional diffusion equations (FDEs) with fractional derivatives $\alpha\in(1,2)$, using the so-called weighted and shifted Gr\"unwald formula, leads to linear systems whose coefficient matrices show a Toeplitz-like structure.
More precisely, in the case of variable coefficients, the related matrix sequences belong to the so-called generalized locally Toeplitz (GLT) class. Conversely, when the given FDE has constant coefficients, using a suitable discretization, we encounter a Toeplitz structure  associated to a  nonnegative function $\mathcal{F}_\alpha$,  called the spectral symbol, having a unique zero at zero of real positive order between one and two. For the fast solution of such systems by preconditioned Krylov methods, several  preconditioning techniques have been proposed in both the one- and two-dimensional cases. In this paper we propose a new preconditioner denoted by $\mathcal{P}_{\mathcal{F}_\alpha}$ which belongs to the $\tau$-algebra and it is  based on the spectral symbol $\mathcal{F}_\alpha$. Comparing with some of the previously proposed preconditioners, we show that although the low band structure preserving preconditioners are more effective in the one-dimensional case, the new  preconditioner performs better in the more challenging multi-dimensional setting.
}

\section{Introduction}
\label{sec:introduction}
Fractional calculus may be considered both an old and modern topic. Old, since it dates back to the letter from L'H{\^o}pital to Leibniz in 1695, and a novel one, since  it has been object of specialized conferences and treatises, for the last 40 years. In recent years considerable interest in fractional calculus has been stimulated by the applications that this calculus finds in numerical analysis and modeling. As an example,  fractional diffusion equations (FDEs) are used to model anomalous diffusion or dispersion,  where a particle plume spreads at a rate inconsistent with the classical Brownian motion model (e.g., see \cite{Meerschaert2004} and the references therein). Such phenomena are ubiquitous in both natural and social sciences. In fact, many complex dynamical systems often contain anomalous diffusion. Fractional kinetic equations are usually an effective method for describing these complex systems, including diffusion type, diffusive convection type and Fokker--Planck type FDEs \cite{Saichev}. Since analytical solutions are rarely available, these kinds of equations are of numerical interest. When the order of fractional derivatives  is $\alpha= 1$, we have the standard diffusion process. With  $0 < \alpha < 1$, we describe  a sub-diffusion process or dispersive, slow diffusion process with the anomalous diffusion index,  while with $\alpha> 1,$   an ultra-diffusion process or increased, fast diffusion process. In \cite{donatelli161}  it has been proved   the strict relationship between  the  order of the fractional derivative and the order of the zero of the associated  symbol of the coefficient matrix of the related system. In addition, it is well known (e.g., see \cite{serraLA:98,serra981}), that  if the generating function $f\ge 0$ of a Toeplitz matrix of size $n$ has a unique zero of order $\alpha\in (0,1)$ and it is bounded, then $T_n(f)$ has a condition number growing exactly as $n^\alpha$.  Hence, the standard Conjugate Gradient method requires only $O(n^{\alpha/2})$ iterations for the solution of a related linear systems up to the required accuracy (for a total of $O(n^{\alpha/2}\log(n))$ arithmetic operations.   When $\omega f$ has all the range in the right complex plane for some $\omega$ complex of modulus one, the generating function $f$ is called weakly sectorial~\cite{bottcher981,sertilli}.  Then if $f$ is weakly sectorial, essentially bounded and $f$ has a unique zero of order $\alpha \in (0,1)$, then again $T_n(f)$ has a condition number growing exactly as $n^\alpha$. Hence, a good GMRES with (possibly)  any standard circulant preconditioning is essentially satisfactory (e.g., see \cite{channg,ng,serhuc}). Thus, the case where the order of fractional derivative $\alpha$ belongs to the interval  $(1,2)$ is, computationally, more challenging.
Moreover, in this paper we are focus on the numerical solution of particular time-dependent space-fractional diffusion equation on rectangular domains in one and two dimensions using finite differences techniques.  For  numerical  techniques  concerning  domains of general geometry or numerical schemes different from finite differences and multigrid  techniques, the interested reader is referred, for example, to  \cite{BKM}, \cite{HLM}, and \cite{KW},  and the references therein. 

Several definitions for the fractional derivative exist, and each  of them approaches the definition of  ordinary derivative in the integer order limit. In \cite{Meerschaert2004,Meerschaert2006} the authors proposed two unconditionally stable finite difference schemes, of first and second order accuracy, based on the shifted Gr{\"u}nwald--Letnikov definition of fractional derivatives. 
In \cite{Wang2010} it was shown that once one of these methods is chosen, the coefficient matrix of the generated system can be seen as the sum of two structures, each of them expressed as a diagonal matrix multiplied by  a Toeplitz matrix.
Since the efficient solution of such systems are of great interest many iterative solvers have been proposed.
Representative   examples are the   multigrid method (MGM) scheme proposed by \cite{Pang2012}, the circulant-based preconditioners  for the  Conjugate Gradient Normal Residual (CGNR) method \cite{lei131,PNW:16}, the splitting preconditioner  \cite{LNS:17},  and two structure-preserving preconditioners proposed in  \cite{donatelli161}.
In the latter paper, the authors provide a detailed analysis, showing that the sequence of coefficient matrices belongs to the generalized locally Toeplitz (GLT) class and its spectral symbol, which describes the asymptotic singular and eigenvalue distribution, is explicitly derived.
In \cite{moghaderi171} the analysis is extended to the two-dimensional case and the authors compare the two-dimensional version of the structure preserving preconditioner based on a decomposition of the Laplacian \cite{donatelli161} to a preconditioner based on an algebraic MGM.

In this work, based on  the theoretical  results presented in \cite{noutsos161} and  motivated by an interest to study  the effectiveness of suitable $\tau$ preconditioners for ill-conditioned symmetric Toeplitz systems, we propose a new preconditioner for the solution of Toeplitz-like systems, stemming from the discretization of the considered FDEs.  Specifically,  in   \cite{noutsos161} the  authors proved the essential spectral equivalence between the matrix sequences $\{T_n(f)\}_n$ and $\{\tau_n(f)\}_n,$  where  $\{T_n(f)\}_n$ is the sequence of symmetric positive definite (SPD) Toeplitz matrices generated by an even, non-negative functions $f$ with zeros  of any positive  order, i.e.,

\begin{displaymath}
[T_n(f)]_{kj}=[T_n(f)]_{k-j}=\frac{1}{2\pi}\int_{-\pi}^{\pi}f(x){\rm e}^{-\boldsymbol{i}(k-j)x}  \;dx \qquad k,j=1,2, \ldots , n, \quad \boldsymbol{i}^2=-1,
\end{displaymath}
and $\{\tau_n(f)\}_n$ is the sequence of a specific  $\tau$ matrices, generated as 
\begin{displaymath}
\tau_n(f)=\mathbb{S}_n\mathrm{diag}(f(\boldsymbol{ \theta}))\mathbb{S}_n, \qquad \boldsymbol{ \theta}=\left[\theta_1, \theta_2,\ldots, \theta_n  \right], \qquad  \theta_j= \frac{j\pi}{n+1}=j\pi h, \qquad j=1,\ldots,n, \label{tau_preconditioner}
\end{displaymath}
and
\begin{align}
[\mathbb{S}_n]_{i,j}&=\sqrt{\frac{2}{n+1}}\sin{\left(i\theta_j\right)},  \qquad i,j=1,\ldots, n.  \label{sine_transform}
\end{align}
We recall here that $\mathbb{S}_n$ is symmetric and orthogonal and so it coincides with its inverse  and that  `essential spectral equivalence' means that all the eigenvalues of $\{\tau_n^{-1}(f)T_n(f)\}_n$ belong to an interval $[c, C]$ except possible $m$ outliers, not converging to zero as the matrix size tends to infinity. In the case of generating functions with the order of their zero lying in the interval $[0,3]$ it is worth noticing that there are no outliers.

According to the analysis given in the aforementioned works, the coefficient matrix of the corresponding linear system depends on  the diffusion coefficients  of the FDE. In the simplest case where they  are constant and equal, the related matrix is a diagonal matrix multiplied  by a real SPD Toeplitz matrix with its generating function $\mathcal{F}_\alpha$ being  even, positive, and real, having  a zero at zero of real positive order between one and two, plus a positive diagonal with constant entries that asymptotically tend to zero.
Analysis shows that this matrix is present in the more general case where the diffusion coefficients are not constant and  not equal to each other. In this case,  a diagonal times skew-symmetric real Toeplitz matrix is then added to the coefficient matrix.
Taking advantage of this fact, we propose  the  preconditioner $\mathcal{P}_{\mathcal{F}_\alpha}=D_n\tau_n(\mathcal{F}_\alpha)$, where $D_n$ is a suitable diagonal matrix defined in Section \ref{sec:main}.
We show that this preconditioner can effectively keep the real part of the eigenvalues away from zero, while the sine transform keeps the cost per iteration $\mathcal{O}(n\log n)$, using a specific real algorithm or using the fast Fourier transform (FFT).
It turns out  that this preconditioner is very efficient and performs better,  especially  in multi-dimensional case, than the proposed preconditioners in \cite{donatelli161} and  \cite{moghaderi171}.

The paper is organized as follows.
In Sections~\ref{sec:introduction:fde}--\ref{sec:introduction:fde2d} we present the one  and two-dimensional FDE problems and the respective discretizations.
Then, in Section~\ref{sec:spectral} we summarize the spectral analysis performed in \cite{donatelli161,moghaderi171},  which turns out to be necessary for the definition of the new preconditioner.
In Section~\ref{sec:main} we also define the proposed preconditioners for the one and two-dimensional cases.
In Section~\ref{sec:numerical} we report numerical experiments and results that confirm the efficiency of the proposed preconditioner.
Finally, in Section~\ref{sec:conclusions} we discuss the advantages and disadvantages of the proposed preconditioners and possible future research directions.

\subsection{Fractional diffusion equations}
\label{sec:introduction:fde}
Consider the two dimensional initial-boundary value problem
\begin{align}
\begin{cases}
\frac{\partial u(x,y,t)}{\partial t}=
d_+(x,y,t)\frac{\partial^\alpha u(x,y,t)}{\partial_+ x^\alpha}+
d_-(x,y,t)\frac{\partial^\alpha u(x,y,t)}{\partial_- x^\alpha}+\\
\hspace{1.6cm}+ e_+(x,y,t)\frac{\partial^\beta u(x,y,t)}{\partial_+ y^\beta}+
e_-(x,y,t)\frac{\partial^\beta u(x,y,t)}{\partial_- y^\beta}+
f(x,y,t),&(x,y,t)\in\Omega\times (0,T),\\
u(x,y,t)=0,&(x,y,t)\in\mathbb{R}^2\setminus\Omega\times[0,T],\\
u(x,y,0)=u_0(x,y),&(x,y)\in\bar{\Omega},
\end{cases}
\label{eq:fde}
\end{align}
where $\Omega=(L_1,R_1)\times (L_2,R_2),$ $\alpha,\beta \in (1,2)$ is the fractional derivative order, $f(x,y,t)$ is the source term and the nonnegative functions $d_{\pm}(x,y,t)$ and $e_{\pm}(x,y,t)$ are the diffusion coefficients. Accordingly, in the one-dimensional setting we drop the dependency on $y$, while the terms including $e_{\pm}(x,y,t)$ are not present.

The left-handed ($\partial_+$) and the right-handed ($\partial_-$) fractional derivatives in \eqref{eq:fde} are defined in Riemann--Liouville form as follows:
\begin{align}
\frac{\partial^\alpha u(x,y,t)}{\partial_+ x^\alpha}&=
\frac{1}{\Gamma(2-\alpha)}\frac{\partial^2}{\partial x^2}\int_{L_1}^x\frac{u(\xi,y,t)}{(x-\xi)^{\alpha-1}}\mathrm{d}\xi,
&\frac{\partial^\alpha u(x,y,t)}{\partial_- x^\alpha}&=
\frac{1}{\Gamma(2-\alpha)}\frac{\partial^2}{\partial x^2}\int_x^{R_1}\frac{u(\xi,y,t)}{(\xi-x)^{\alpha-1}}\mathrm{d}\xi,\nonumber\\
\frac{\partial^\beta u(x,y,t)}{\partial_+ y^\beta}&=
\frac{1}{\Gamma(2-\beta)}\frac{\partial^2}{\partial y^2}\int_{L_2}^y\frac{u(x,\eta,t)}{(y-\eta)^{\beta-1}}\mathrm{d}\eta,
&\frac{\partial^\beta u(x,y,t)}{\partial_- y^\beta}&=
\frac{1}{\Gamma(2-\beta)}\frac{\partial^2}{\partial y^2}\int_y^{R_2}\frac{u(x,\eta,t)}{(\eta-y)^{\beta-1}}\mathrm{d}\eta.\nonumber
\end{align}

\subsection{First-order finite difference discretization}
\label{sec:introduction:fde1}
In this section, we consider the one-dimensional version of~\eqref{eq:fde} (for two-dimensional derivation see Section~\ref{sec:introduction:fde2d} and~\cite{moghaderi171}). Applying the shifted Gr\"unwald formulas we can approximate the left and right fractional derivatives by
\begin{align}
\frac{\partial^\alpha u(x,t)}{\partial_+ x^\alpha}&=\frac{1}{h_x^\alpha}\sum_{k=0}^{\lfloor(x-L_1)/h_x\rfloor}g_k^{(\alpha)}u(x-(k-1)h_x,t)+\mathcal{O}(h_x),\nonumber\\
\frac{\partial^\alpha u(x,t)}{\partial_- x^\alpha}&=\frac{1}{h_x^\alpha}\sum_{k=0}^{\lfloor(R_1-x)/h_x\rfloor}g_k^{(\alpha)}u(x+(k-1)h_x,t)+\mathcal{O}(h_x),\nonumber
\end{align}
where $\lfloor \cdot\rfloor$ is the floor function, $n_1$ is the discretization parameter giving $h_x=(R_1-L_1)/(n_1+1)=(R_1-L_1)h_1$,  and $g_k^{(\alpha)}$ are the alternating fractional binomial coefficients defined as
\begin{align}
g_k^{(\alpha)}=(-1)^k\binom{\alpha}{k}=\frac{(-1)^k}{k!}\alpha(\alpha-1)\cdots(\alpha-k+1),\quad k= 0,1,\ldots,\label{eq:grunwald1}
\end{align}
where $\binom{\alpha}{0}=1$.
Using the implicit Euler method for time discretization, we define the number of time steps (index $m$) to be $M$, and thus $h_t=T/M$, and 
\begin{align}
\frac{u_i^{(m)}-u_i^{(m-1)}}{h_t}=\frac{d_{+,i}^{(m)}}{h_x^\alpha}\sum_{k=0}^{i+1}g_k^{(\alpha)}u_{i-k+1}^{(m)}+\frac{d_{-,i}^{(m)}}{h_x^\alpha}\sum_{k=0}^{n_i-i+2}g_k^{(\alpha)}u_{i+k-1}^{(m)}+f_i^{(m)},\nonumber
\end{align}
where $d_{\pm,i}^{(m)}=d_\pm(x_i,t_m)$, $u_i^{(m)}=u(x_i,t_m)$, and $f_i^{(m)}=f(x_i,t_m)$, where $x_i=L_1+ih_x$ and $t_m=mh_t$.
After rearranging terms, we find
\begin{align}
\frac{h_x^\alpha}{h_t}u_i^{(m)}-d_{+,i}^{(m)}\sum_{k=0}^{i+1}g_k^{(\alpha)}u_{i-k+1}^{(m)}-d_{-,i}^{(m)}\sum_{k=0}^{n_1-i+2}g_k^{(\alpha)}u_{i+k-1}^{(m)}=\frac{h_x^\alpha}{h_t}u_i^{(m-1)}+h_x^\alpha f_i^{(m)}, \nonumber
\end{align}
or in matrix form, the linear systems
\begin{align}
\left(\nu_{M,n_1}\mathbb{I}_{n_1}+D_+^{(m)}T_{\alpha,n_1}+D_-^{(m)}T_{\alpha,n_1}^{\mathrm{T}}\right)\mathbf{u}^{(m)}&=\nu_{M,n_1}\mathbf{u}^{(m-1)}+h_x^\alpha\mathbf{f}^{(m)},
\label{eg:matrixform}
\end{align}
where
\begin{align}
\mathbb{I}_{n_1}&: \mbox{The identity  matrix of size} ~n_1,\\ 
\nu_{M,n_1}&=\frac{h_x^\alpha}{h_t}, \label{nu_definition}\\
\mathbf{u}^{(m)}&=\left[u_1^{(m)},u_2^{(m)},\ldots,u_{n_1}^{(m)}\right]^{\mathrm{T}},\nonumber\\
\mathbf{f}^{(m)}&=\left[f_1^{(m)},f_2^{(m)},\ldots,f_{n_1}^{(m)}\right]^{\mathrm{T}},\nonumber\\
[D_{\pm}^{(m)}]_{i,i}&=d_{\pm}^{(m)}(x_i,t_m),\quad i=1,\ldots,n_1,\nonumber
\end{align}
and
\begin{align}
T_{\alpha,n_1}=-\left[
\begin{array}{ccccccccc}
g_1^{(\alpha)}&g_0^{(\alpha)}&\phantom{\ddots}\\
g_2^{(\alpha)}&g_1^{(\alpha)}&g_0^{(\alpha)}&\phantom{\ddots}\\
g_3^{(\alpha)}&g_2^{(\alpha)}&g_1^{(\alpha)}&g_0^{(\alpha)}&\phantom{\ddots}\\
\vdots&\ddots&\ddots&\ddots&\ddots\\
\vdots&\ddots&\ddots&\ddots&\ddots&\ddots\\
g_{n_1-1}^{(\alpha)}&g_{n_1-2}^{(\alpha)}&\cdots&\ddots&g_2^{(\alpha)}&g_1^{(\alpha)}&g_0^{(\alpha)}\\
g_{n_1}^{(\alpha)}&g_{n_1-1}^{(\alpha)}&\cdots\vphantom{\ddots}&\cdots&g_3^{(\alpha)}&g_2^{(\alpha)}&g_1^{(\alpha)}\\
\end{array}
\right],
\label{eq:grunwaldmatrix1}
\end{align}
with the coefficients $g_k^{(\alpha)}$  given in \eqref{eq:grunwald1}.

Now define
\begin{align}
\mathcal{M}_{\alpha,n_1}^{(m)}&=\left(\nu_{M,n_1}\mathbb{I}_{n_1}+D_+^{(m)}T_{\alpha,n_1}+D_-^{(m)}T_{\alpha,n_1}^{\mathrm{T}}\right),\label{first_order_1d_coefficient_matrix} \\
\mathbf{b}^{(m)}&=\mathbf\nu_{M,n_1}{u}^{(m-1)}+h_x^\alpha\mathbf{f}^{(m)}. \nonumber
\end{align}
Then, for each time step $m,$ we solve the system
\begin{align}
\mathcal{M}_{\alpha,n_1}^{(m)}\mathbf{u}^{(m)}&=\mathbf{b}^{(m)}.\label{eq:linsystem}
\end{align}

\subsection{Second-order finite difference discretization}
\label{sec:introduction:fde2}
For the second order finite difference discretization in space, we can just exchange the matrix $T_{\alpha,n_1}$ in \eqref{eg:matrixform} with a matrix $S_{\alpha,n_1}$ defined by
\begin{align}
S_{\alpha,n_1}=-\left[
\begin{array}{ccccccccc}
w_1^{(\alpha)}&w_0^{(\alpha)}&\phantom{\ddots}\\
w_2^{(\alpha)}&w_1^{(\alpha)}&w_0^{(\alpha)}&\phantom{\ddots}\\
w_3^{(\alpha)}&w_2^{(\alpha)}&w_1^{(\alpha)}&w_0^{(\alpha)}&\phantom{\ddots}\\
\vdots&\ddots&\ddots&\ddots&\ddots\\
\vdots&\ddots&\ddots&\ddots&\ddots&\ddots\\
w_{n_1-1}^{(\alpha)}&w_{n_1-2}^{(\alpha)}&\cdots&\ddots&w_2^{(\alpha)}&w_1^{(\alpha)}&w_0^{(\alpha)}\\
w_{n_1}^{(\alpha)}&w_{n_1-1}^{(\alpha)}&\cdots\vphantom{\ddots}&\cdots&w_3^{(\alpha)}&w_2^{(\alpha)}&w_1^{(\alpha)}\\
\end{array}
\right],
\label{eq:grunwaldmatrix2}
\end{align}
where
\begin{align}
    w_0^{(\alpha)}&=\frac{\alpha}{2}g_0^{(\alpha)},\nonumber\\
    w_k^{(\alpha)}&=\frac{\alpha}{2}g_k^{(\alpha)}+\frac{2-\alpha}{2}g_{k-1}^{(\alpha)}, \quad k\geq 1, \nonumber
\end{align}
and the coefficients $g_k^{(\alpha)}$ are expressed as in relation \eqref{eq:grunwald1}.

\subsection{Two-dimensional case}
\label{sec:introduction:fde2d}
Similarly  to 1D case, we can extend the discretization scheme to the  two-dimensional setting. In the next paragraph we  summarize the main points  of the numerical procedure, referring the reader   in \cite{moghaderi171} for further details. Define
\begin{align}
h_x&=\frac{R_1-L_1}{n_1+1}=(R_1-L_1)h_1,\quad\quad x_{i}=L_1+ih_x,\quad i=1,\ldots, n_1,\nonumber\\
h_y&=\frac{R_2-L_2}{n_2+1}=(R_2-L_2)h_2,\quad\quad y_{i}=L_2+ih_y,\quad i=1,\ldots, n_2,\nonumber
\end{align}
and $N=n_1n_2$.
The solution $u(x,y,t)$ is discretized as $u_{i,j}^{(m)}=u(x_i,y_j,t^{(m)})$,
\begin{align}
\mathbf{u}^{(m)}&=[u_{1,1}^{(m)},\ldots,u_{n_1,1}^{(m)},u_{1,2}^{(m)},\ldots,u_{n_1,2}^{(m)},\ldots, u_{1,n_2}^{(m)},\ldots,u_{n_1,n_2}^{(m)}]^{\mathrm{T}},\nonumber
\end{align}
and the four diffusion function $d_+(x,y,t)$, $d_-(x,y,t)$,  $e_+(x,y,t)$, $e_-(x,y,t)$ are discretized as $d_{i,j}^{\pm,(m)}=d_\pm(x_i,y_j,t^{(m)})$ and $e_{i,j}^{\pm,(m)}=e_\pm(x_i,y_j,t^{(m)}),$
\begin{align}
\mathbf{d}_{\pm}^{(m)}&=[d_{1,1}^{\pm,(m)},\ldots,d_{n_1,1}^{\pm,(m)},d_{1,2}^{\pm,(m)},\ldots,d_{n_1,2}^{\pm,(m)},\ldots, d_{1,n_2}^{\pm,(m)},\ldots,d_{n_1,n_2}^{\pm,(m)}]^{\mathrm{T}},\nonumber\\
\mathbf{e}_{\pm}^{(m)}&=[e_{1,1}^{\pm,(m)},\ldots,e_{n_1,1}^{\pm,(m)},e_{1,2}^{\pm,(m)},\ldots,e_{n_1,2}^{\pm,(m)},\ldots, e_{1,n_2}^{\pm,(m)},\ldots,e_{n_1,n_2}^{\pm,(m)}]^{\mathrm{T}}.\nonumber
\end{align}
The source term $f(x,y,t)$ is discretized as $f_{i,j}^{(m)}=f(x_i,y_j,t^{(m)})$,
\begin{align}
\mathbf{v}^{(m-1/2)}&=[f_{1,1}^{(m-1/2)},\ldots,f_{n_1,1}^{(m-1/2)},f_{1,2}^{(m-1/2)},\ldots,f_{n_1,2}^{(m-1/2)},\ldots, f_{1,n_2}^{(m-1/2)},\ldots,f_{n_1,n_2}^{(m-1/2)}]^{\mathrm{T}}.\nonumber
\end{align}
We also  define the four matrices $D_\pm^{(m)}=\mathrm{diag}(\mathbf{d}_\pm^{(m)})$ and $E_\pm^{(m)}=\mathrm{diag}(\mathbf{e}_\pm^{(m)})$.

If we have  two fractional derivatives, $\alpha$ and $\beta$, in each spatial direction  we define the two matrices $S_{\alpha,n_1}$ and $S_{\beta,n_2}$ (or $T_{\alpha,n_1}$ and $T_{\beta,n_2}$ for the considered  first-order discretization).

We also define the two $N\times N$ matrices
\begin{align}
A_x^{(m)}&=D_+^{(m)}(\mathbb{I}_{n_2}\otimes S_{\alpha,n_1})+D_-^{(m)}(\mathbb{I}_{n_2}\otimes S_{\alpha,n_1}^{\mathrm{T}}),\nonumber\\
A_y^{(m)}&=E_+^{(m)}( S_{\beta,n_2}\otimes \mathbb{I}_{n_1})+E_-^{(m)}( S_{\beta,n_2}^{\mathrm{T}}\otimes \mathbb{I}_{n_1}),\nonumber
\end{align}
where $\mathbb{I}_n$ denotes the identity matrix of size $n$, and $\otimes$ is the Kronecker product. Using Crank--Nicolson approach for time discretization (e.g., see \cite{moghaderi171}) we obtain the system

\begin{align}
    \left(\frac{1}{r}\mathbb{I}_{N}+A_{x}^{(m)}+\frac{s}{r}A_{y}^{(m)}\right)\mathbf{u}^{(m)}=\left(\frac{1}{r}\mathbb{I}_{N}-A_{x}^{(m-1)}-\frac{s}{r}A_{y}^{(m-1)}\right)\mathbf{u}^{(m-1)}+2h_x^{\alpha}\mathbf{v}^{(m-1/2)}\nonumber
\end{align}
where $r=\frac{h_t}{2h_x^\alpha}$, $s=\frac{h_t}{2h_y^\beta}$. In compact form we have
\begin{align}
\mathcal{M}_{(\alpha,\beta),N}^{(m)}\mathbf{u}^{(m)}=\mathbf{b}^{(m)},\nonumber
\end{align}
where
\begin{align}
\mathcal{M}_{(\alpha,\beta),N}^{(m)}&=\frac{1}{r}\mathbb{I}_{N}+A_{x}^{(m)}+\frac{s}{r}A_{y}^{(m)},\nonumber\\
\mathbf{b}^{(m)}&=\left(\frac{1}{r}\mathbb{I}_{N}-A_{x}^{(m-1)}-\frac{s}{r}A_{y}^{(m-1)}\right)\mathbf{u}^{(m-1)}+2h_x^{\alpha}\mathbf{v}^{(m-1/2)}.\nonumber
\end{align}

\section{Spectral analysis}
 
In this section we provide  some definitions that are used in the analysis. We also employ the theory of GLT matrix sequences to study the spectral properties of $\mathcal{M}_{\alpha,n_1}^{(m)}$ of \eqref{eq:linsystem} (for both the first and second order version) as the matrix dimension tends to infinity. We refer the reader to \cite{garoni171} for an introduction to the theory of GLT matrix sequences. Here, we only list some basic properties that are used in the analysis that follows. The results reported in Section \ref{ssec:sp1} and Section \ref{ssec:sp2} are taken from \cite{donatelli161,moghaderi171}.

\label{sec:spectral}
\begin{definition}
Let $\{A_n\}_n$ be a matrix sequence and $f:D\rightarrow \mathbb{C}$ be a measurable function defined on a measurable set $D\subset \mathbb{R}^k$ with $0<\mu(D)<\infty$. 
\begin{itemize}
    \item We say that the sequence $\{A_n\}_n$ has an asymptotic singular value distribution described by $f$, and we write $\{A_n\}_n \sim_{\sigma}f$ if,  
    \[
    \lim_{n\to \infty} \frac{1}{n} \sum_{j=1}^{n}F(\sigma_j(A_n)) = \frac{1}{\mu(D)}\int_DF(|f(x)|)dx, \quad \forall F \in C_c(\mathbb{R}),
    \]
    where $C_c(\mathbb{R})$ is the set of continuous functions with compact support over $\mathbb{R}$.

    \item We say that $\{A_n\}_n$ has an asymptotic eigenvalue  distribution described by $f$, and  write $\{A_n\}_n \sim_{\lambda}f$ if
    \[
    \lim_{n\to \infty} \frac{1}{n} \sum_{j=1}^{n}F(\lambda_j(A_n)) = \frac{1}{\mu(D)}\int_DF(f(x))dx, \quad \forall F \in C_c(\mathbb{C}),
    \]
    where $C_c(\mathbb{R})$ is the set of continuous functions with compact support over $\mathbb{R}$.
\end{itemize}
\end{definition}

\begin{definition}
Let $f\in L^{1}([-\pi,\pi])$ and $\{f_k\}_{k\in\mathbb{Z}}$ its Fourier coefficients defined as
\[
  f_{k}=\frac{1}{2\pi} \int_{-\pi} ^{\pi}f(\theta) e^{-\mathbf{i} k \theta } \,d\theta,\quad k=0,\pm1,\pm2,\dots.
\]
The sequence of matrices $\{T_n(f)\}_{n\in\mathbb{N}}$, $T_n(f)=[f_{i-j}]_{i,j=1}^{n}$, is called  a Toeplitz sequence generated by $f$. 
\end{definition}
The eigenvalue and singular value distribution of Toeplitz sequences generated by $f\in L^1([-\pi, \pi])$ is given by generalized Szeg{\H{o}} theorem \cite{MR1481397}:

  \begin{theorem}\label{thm:toeplitz spectral symbol}
    Let $f \in L^{1}([-\pi,\pi])$ and $T_n(f)$ be the Toeplitz matrix generated by $f$. Then $f$ is the spectral symbol of the sequence, that is 
    \[
    \{T_n(f)\}_n \sim_{\sigma} f.
    \]
    If, moreover, $f$ is real-valued, then
    \[
    \{T_n(f)\}_n \sim_{\lambda} f.
    \]
  \end{theorem}

The basic properties of the GLT class follow.
\begin{itemize}
   \item[GLT1] Each GLT sequence $\{A_n\}_n$ has a singular value symbol $\tilde f:[0,1]\times [-\pi,\pi] \rightarrow \mathbb{C}$. If all the matrices of the sequence are Hermitian, then the distribution also holds in the eigenvalue sense. We call $\tilde f(x,\theta)$ the GLT symbol of $\{A_n\}_n$ and we write $\{A_n\}_n \sim_{\textsc{glt}} \tilde f$.
   \item[GLT2]  The set of GLT sequences  is closed under linear combinations, products,  inversion (whenever the symbol is zero in at most a set of zero Lebesgue measure) and conjugation. The sequence  obtained via algebraic operations on a finite set of given GLT sequences is still a GLT sequence and its symbol is  obtained by performing the same algebraic manipulations on the corresponding symbols of the input GLT sequences.
   \item[GLT3] Every Toeplitz sequence generated by a function $f \in L^1([-\pi,\pi])$ is a GLT sequence and its symbol is $\tilde{f}(x,\theta) = f(\theta)$. If $\alpha:[0,\:1]\rightarrow\mathbb{C}$ is a Riemman integrable function, the diagonal matrix sequence of the form $\{D_n(\alpha)\}_n, \: n\in \mathbb{N}$,  $D_n(\alpha)=\mbox{diag}_{j=1,\dots,n}(\alpha(\frac{j}{n}))$ is a GLT sequence with spectral symbol $\tilde{f}(x,\theta) = \alpha(x)$.
 \end{itemize}

\subsection{Spectral analysis: Matrices $T_{\alpha,n}$ and $S_{\alpha,n}$}\label{ssec:sp1}

From~\cite{donatelli161} we know that  $T_{\alpha,n_1}$ in \eqref{eq:grunwaldmatrix1} is a Toeplitz sequence with spectral symbol
\begin{align}
g_\alpha(\theta)=-e^{-\mathbf{i}\theta}\left(1-e^{\mathbf{i}\theta}\right)^\alpha,\nonumber
\end{align}
and thus from Theorem \ref{thm:toeplitz spectral symbol} and GLT3
\begin{align}
    \{T_{\alpha,n}\}_{n}=\{T_{n}(g_\alpha)\}_{n}\sim_{\textsc{glt},\sigma}g_\alpha. \label{T_as_Toeplitz}
\end{align}
Furthermore, as shown in~\cite{moghaderi171}, $S_{\alpha,n}$ in \eqref{eq:grunwaldmatrix2} is a Toeplitz sequence with spectral symbol
\begin{align}
w_\alpha(\theta)=-\left(\frac{2-\alpha(1-e^{-\mathbf{i}\theta})}{2}\right)\left(1-e^{\mathbf{i}\theta}\right)^\alpha,\label{eq_w_function}
\end{align}
and thus from Theorem \ref{thm:toeplitz spectral symbol} and GLT3
\begin{align}
   \{S_{\alpha,n}\}_n= \{T_n(w_\alpha)\}_n\sim_{\textsc{glt},\sigma}w_\alpha. \label{S_as_Toeplitz}
\end{align}

\subsection{Spectral analysis: Constant coefficient case}\label{ssec:sp2}
\begin{theorem}
\label{thm:const}
Assuming $d_{\pm}(x,t)=d>0$ and that $\nu_{M,n}=o(1)$,  we have for the first order spatial discretization
\begin{align}
\{\mathcal{M}_{\alpha,n}^{(m)}\}_n\sim_{\textsc{glt},\sigma,\lambda}d\cdot p_\alpha(\theta),\nonumber
\end{align}
where
\begin{align}
p_{\alpha}(\theta)&=g_\alpha(\theta)+g_\alpha(-\theta).\label{first_order_symbol}
\end{align}
For the second order spatial discretization we have
\begin{align}
\{\mathcal{M}_{\alpha,n}^{(m)}\}_n\sim_{\textsc{glt},\sigma,\lambda}d\cdot q_\alpha(\theta),\nonumber
\end{align}
where
\begin{align}
q_{\alpha}(\theta)&=w_\alpha(\theta)+w_\alpha(-\theta).\label{second_order_symbol}
\end{align}
In the subsequent analysis,  whenever either symbol  $p_\alpha(\theta)$ or $q_\alpha(\theta)$ is applicable, we denote both symbols by $\mathcal{F}_\alpha(\theta)$.
\end{theorem}
\begin{proposition}\label{prop:sp2}
Let $\alpha \in (1,2)$, then the function $p_{\alpha}(\theta)$ has a zero of order $\alpha$ at 0.
\end{proposition}
Moreover, in connection with Proposition \ref{prop:sp2}, it is worth noticing the following: if $f$ is nonnegative with a unique zero of order $\alpha>0$, then the matrix
$T_n(f)$ is positive definite for any $n$ its minimal eigenvalue tends to zero as $n$ tends to infinity as $n^{-\alpha}$; furthermore, if $f$ is also bounded then the condition number
of $T_n(f)$ grows asymptotically as $n^{\alpha}$ (e.g., see \cite{bottcher981, serra981}).

\section{Main Results}
\label{sec:main}
In this section we propose two new preconditioners, based on the spectral symbol, for the one and two-dimensional problems.
\subsection{Proposed Preconditioner: One Dimension}
\label{sec:prop1d}
To be consistent with \cite{donatelli161}, so that results can be compared, we use the first-order spatial discretization for the one dimensional case.
We   also omit the time dependency mark to simplify the notation. Thus, let $T_n=T_{\alpha,n1}$ be defined as in \eqref{eq:grunwaldmatrix1} and let $\mathcal{M}_{n}=\mathcal{M}_{\alpha,n_1}$ be defined as in \eqref{first_order_1d_coefficient_matrix}.

As previously mentioned in Section~\ref{sec:introduction}, the proposed preconditioner is similar to  a diagonal matrix $D_n$ times a specific $\tau$ matrix, i.e., $\mathcal{P}_{\mathcal{F}_\alpha}=D_n\tau_n(\mathcal{F}_{\alpha}(\boldsymbol{\theta}))$, where both these  parts  will be clarified through this paragraph.

  The product of two or more matrices  as preconditioner is not a new proposal (see, e.g., \cite{NV-2008}).
The coefficient matrix of the system $\mathcal{M}_n=\nu_{M,n}\mathbb{I}_{n}+D_+T_{n}+D_-T_{n}^{\mathrm{T}}$ suggests the following candidate for the diagonal matrix
\begin{align}
D_n&=\frac{1}{2}\left(D_++D_-\right), \nonumber \\
[D_n]_{i,i}&= \frac{d_{+,i}+d_{-,i}}{2}, \label{Diagonal_1D_matrix}
\end{align}
that has been used in other preconditioning strategies (see  for example \cite{donatelli161}). Then, assuming that $d_{\pm}$ do not have a common zero at  $x_0 \in [L,R]$, we deduce that $D_n^{-1}$ is uniformly bounded and
\begin{align}
D_n^{-1}\mathcal{M}_n=\nu_{M,n}D_n^{-1}+D_n^{-1}D_+T_{n}+D_n^{-1}D_-T_{n}^{\mathrm{T}}.\nonumber
\end{align}
Defining $\delta(x)=\frac{d_+(x)}{d_+(x)+d_-(x)}$, $\delta_i=\delta(x_i)$, $\boldsymbol{\delta}=[\delta_1,\delta_2,\ldots,\delta_n]$, and $G_n=\mathrm{diag}(\boldsymbol{\delta}),$ taking into account that $d_\pm$ are non negative  functions, we have that $0\leq \delta(x)\leq 1$ and also
\begin{align}
D_n^{-1}D_+&=2G_n, \nonumber\\
D_n^{-1}D_-&=2(\mathbb{I}_{n}-G_n).\nonumber
\end{align}
Hence, $D_n^{-1}\mathcal{M}_n$ can be written as
\begin{align}
D_n^{-1}\mathcal{M}_n&=\nu_{M,n}D_n^{-1}+D_n^{-1}D_+T_{n}+D_n^{-1}D_-T_{n}^{\mathrm{T}}\nonumber\\
&=\nu_{M,n}D_n^{-1}+2G_nT_{n}+2(\mathbb{I}_{n}-G_n)T_{n}^{\mathrm{T}} \nonumber \\
&=\nu_{M,n}D_n^{-1}+(T_{n}+T_{n}^{\mathrm{T}})+(2G_n-\mathbb{I}_{n})(T_n-T_{n}^{\mathrm{T}}). \nonumber
\end{align}
Since, from \eqref{T_as_Toeplitz}, $T_n\coloneqq T_n(-e^{-\mathbf{i}\theta}\left(1-e^{\mathbf{i}\theta}\right)^\alpha)=T_n(g_{\alpha}(\theta))$ and  $T_n^{\mathrm{T}}\coloneqq T_n(-e^{\mathbf{i}\theta}\left(1-e^{-\mathbf{i}\theta}\right)^\alpha)=T_n(g_{\alpha}(-\theta))$ we have
\begin{align}
D_n^{-1}\mathcal{M}_n&=\nu_{M,n}D_n^{-1}+(T_{n}+T_{n}^{\mathrm{T}})+(2G_n-\mathbb{I}_{n})(T_n-T_{n}^{\mathrm{T}}) \nonumber \\
&=\nu_{M,n}D_n^{-1}+T_n(g_{\alpha}(\theta)+g_{\alpha}(-\theta))+(2G_n-\mathbb{I}_{n})T_n(g_{\alpha}(\theta)-g_{\alpha}(-\theta)) \nonumber  \\
&=\nu_{M,n}D_n^{-1}+T_n(p_{\alpha}(\theta))+(2G_n-\mathbb{I}_{n})T_n(2\mathbf{i}\Im\left\{g_{\alpha}(\theta)\right\}),   \label{DM_analytical_presentation}
\end{align}
where  $p_{\alpha}(\theta)$, defined in \eqref{first_order_symbol}, is real, positive and even. With $\Im$ we denote the imaginary part of a function. 
The above derivation of the $D_n^{-1}\mathcal{M}_n$ matrix is of interest since it makes clear why it is reasonable to use the $\tau$ preconditioner.
The first term of the above matrix, $\nu_{M,n}D_n^{-1}$, is  diagonal with positive and $o(1)$ entries, since we have supposed that the $d_{\pm}$ functions do not have zeros at the same point in the domain $[L,R]$ and $\nu_{M,n}=o(1)$.
We mention here that although the entries of this term are $o(1)$, its effect on the eigenvalues of the preconditioned matrix can be  significant. The reason is  explained in the end of this section.
The third term in \eqref{DM_analytical_presentation}  is a diagonal matrix with entries in $[-1,1]$ times a skew-symmetric Toeplitz matrix with generating function $2\mathbf{i}\Im\left\{g_{\alpha}(\theta)\right\}$.
If $d_+=d_-$ this term is vanishing  while if the $d_{\pm}$ are constant but not equal it is a pure skew-symmetric Toeplitz (in that case $(2G_n-\mathbb{I}_{n})=c\mathbb{I}_{n}$ for some constant $c$).

The term in \eqref{DM_analytical_presentation} which is mainly responsible for the dispersion of the real part of the spectrum, is the second term, that is, $T_n(p_{\alpha}(\theta))$. The  $\tau$ preconditioner will effectively  cluster the eigenvalues of this matrix, and consequently the eigenvalues of the whole matrix $D_n^{-1}\mathcal{M}_n$.
Hence, taking advantage of the essential spectral equivalence between the matrix sequences $\{\tau_n(f)\}_n$ and $\{T_n(f)\}_n$ proven in \cite{noutsos161}, we propose a preconditioner expressed as
\begin{align}
\mathcal{P}_{\mathcal{F}_\alpha,n}&=D_n\tau_n(p_{\alpha}(\theta))=D_n\mathbb{S}_nF_n\mathbb{S}_n,  \label{eq:prop_1D}
\end{align}
where

\begin{align}
F_n&=\mathrm{diag}(p_{\alpha}(\boldsymbol{\theta})),\qquad  \boldsymbol{ \theta}=\left[\theta_1, \theta_2,\ldots, \theta_n  \right], \qquad  \theta_j= \frac{j\pi}{n+1}=j\pi h, \qquad j=1,\ldots,n, \nonumber
\end{align}
with $D_n$ defined in \eqref{Diagonal_1D_matrix} and $\mathbb{S}_n$ being the sine transform matrix reported in \eqref{sine_transform}.  

\subsubsection{Case I: $d_{\pm}$ are constants }
  In the case where the diffusion coefficient functions are constants, the \eqref{DM_analytical_presentation} becomes:
  \begin{align*}
  \left(2\frac{\nu_{M,n}}{d_++d_-}\right)\mathbb{I}_n+T_n\left(p_{\alpha}(\theta)\right)+\left(\frac{d_+-d_-}{d_++d_-}\right)T_n\left(2\mathbf{i}\Im\left\{g_{\alpha}(\theta)\right\}\right)=T_n\left(2\frac{\nu_{M,n}}{d_++d_-}+ p_{\alpha}(\theta)\right)+T_n\left(2\left(\frac{d_+-d_-}{d_++d_-}\right)\mathbf{i}\Im\left\{g_{\alpha}(\theta)\right\}\right),
  \end{align*}
  i.e,  is exactly the sum of a symmetric and a skew-symmetric Toeplitz matrix.  
  It is worth noticing that according to the GLT machinery, the term $\frac{2\cdot \nu_{M,n}}{d_++d_-}$ which is added to the symbol  of the first Toeplitz matrix sequence does not change the symbol of the sequence since is of order $o(1)$. However it  affects the speed  in which the minimum eigenvalue of the sequence approaches zero as the dimension of the matrix tends to infinity. Thus, in this special case,  the $\tau$ part of preconditioner  is defined as 
  \begin{align*}
  \tau_{M,n}\left(p_{\alpha}(\theta)+\frac{2\cdot \nu_{M,n}}{d_++d_-}\right)=\mathbb{S}_n \mathrm{diag}\left(p_{\alpha}(\theta)+\frac{2\cdot \nu_{M,n}}{d_++d_-}\right) \mathbb{S}_n=\mathbb{S}_n \hat{F}_n \mathbb{S}_n.
  \end{align*}
  Then,
  \begin{align*}
  \tau_{M,n}^{-1}\left(p_{\alpha}(\theta)+\frac{2\cdot \nu_{M,n}}{d_++d_-}\right) \left[T_n\left(\frac{2\cdot\nu_{M,n}}{d_++d_-}+ p_{\alpha}(\theta)\right)+T_n\left(2\frac{d_+-d_-}{d_++d_-}\mathbf{i}\Im\left\{g_{\alpha}(\theta)\right\}\right)\right] \sim \\
  \hat{F}_n^{-\frac{1}{2}}\mathbb{S}_n\left[ T_n\left(\frac{2\cdot\nu_{M,n}}{d_++d_-}+ p_{\alpha}(\theta)\right)+T_n\left(2\frac{d_+-d_-}{d_++d_-}\mathbf{i}\Im\left\{g_{\alpha}\theta)\right\} \right) \right]\mathbb{S}_n\hat{F}_n^{-\frac{1}{2}}=\\
  \hat{F}_n^{-\frac{1}{2}}\mathbb{S}_nT_n\left(\frac{2\cdot\nu_{M,n}}{d_++d_-}+ p_{\alpha}(\theta)\right)\mathbb{S}_n\hat{F}_n^{-\frac{1}{2}} + \hat{F}_n^{-\frac{1}{2}}\mathbb{S}_nT_n\left(2\frac{d_+-d_-}{d_++d_-}\mathbf{i}\Im\left\{g_{\alpha}(\theta)\right\}\right)\mathbb{S}_n\hat{F}_n^{-\frac{1}{2}}.
  \end{align*}
  The first term in the above sum  is symmetric and its eigenvalues are strongly clustered at 1 since  the conditions of  the main theoretical result of  \cite{noutsos161} are fulfilled concerning   the   spectral equivalence between  a  $\tau$ matrix and  a Toeplitz one. The second term is skew-symmetric and it does not  affect  the real part of the eigenvalues of the whole matrix. Moreover,  it is absent whenever  $d_+=d_-$. Hence, the real parts of the eigenvalues of the preconditioned matrix  are strongly  clustered at  1 and are bounded by constants $c, C$ with $0<c\leq 1\leq C <\infty$.

\subsubsection{Case II.  $d_{-}(x)=d_{+}(x) >0$  }
\label{sec:diffequal}
 In this case,  the term  $2G_n-\mathbb{I}_{n}=\mathbf{0}$ in  \ref{DM_analytical_presentation} is equal to zero and the preconditioned matrix becomes $\tau_n^{-1}(p_{\alpha}(\theta))(\nu_{M,n}D_n^{-1}+T_n(p_{\alpha}(\theta)))$ which is similar to the SPD
\begin{align}
\tau_n^{-1}(p_{\alpha}(\theta))(\nu_{M,n}D_n^{-1}+T_n(p_{\alpha}(\theta)))&\sim F_n^{-1/2}\mathbb{S}_n(\nu_{M,n}D_n^{-1}+T_n(p_{\alpha}(\theta)))\mathbb{S}_nF_n^{-1/2}\nonumber\\
&=\nu_{M,n}F_n^{-1/2}\mathbb{S}_n(D_n^{-1})\mathbb{S}_nF_n^{-1/2}+F_n^{-1/2}\mathbb{S}_n(T_n(p_{\alpha}(\theta)))\mathbb{S}_nF_n^{-1/2}. \label{eq:main:precond1d}
\end{align}
In the above splitting in positive symmetric terms,   the first  one has ${\it o}(n)$ eigenvalues  tending  to infinity while the second one fulfills the main theoretical result of \cite{noutsos161}  and thus, for every $n$, it has eigenvalues belonging to an interval  [c,C] with $c,C$ constants and  $0<c\leq 1\leq C <\infty$. The claim about the spectrum of the first term can be   proved if we equivalently show that the inverse  of it, i.e. $ F_n (\mathbb{S}_n D_n  \mathbb{S}_n)$ has at most ${\it o}(n)$ eigenvalues tending to 0 as $n\rightarrow \infty$. Since $F_n$ is the diagonal matrix formed by the   values  $p_{\alpha}(j\pi h), \quad j=1,\ldots, n,$ which has a unique zero at zero of order $\alpha$,   there will be an index $\hat{j}$ with $\hat{j}$ of order ${\it o}(n)$ such that $p_{\alpha}(j\pi h)$  being  of order ${\it o}(1)$  for all $j\leq \hat{j}$.  Thus, at most ${\it o}(n)$  eigenvalues of $F_n $  can tend to zero. Using  Rayleigh quotient and taking into account that the matrix $D_n$ is a diagonal matrix with entries bounded from above end below by positive universal constants, our claim is proved.   Consequently, using the Weyl's theorem on (\ref{eq:main:precond1d}) we obtain that  
$$
\lambda_{k}\left(\nu_{M,n}F_n^{-1}\mathbb{S}_n(D_n^{-1})\mathbb{S}_n+F_n^{-1}\mathbb{S}_n(T_n(p_{\alpha}(\theta)))\mathbb{S}_n\right) \leq  \nu_{M,n}\lambda_{k}( F_n^{-1}\mathbb{S}_n(D_n^{-1})\mathbb{S}_n) +\lambda_n \left(F_n^{-1}\mathbb{S}_n(T_n(p_{\alpha}(\theta)))\mathbb{S}_n\right). 
$$

Accordingly, at most $\it o(n)$  eigenvalues of $\tau_n^{-1}(p_{\alpha}(\theta))(\nu_{M,n}D_n^{-1}+T_n(p_{\alpha}(\theta)))$ can tend to infinity. Clearly the term $\nu_{M,n}$ which in general tends to zero as ${\it O}(n^{1-\alpha})$, can further  reduce the number of eigenvalues tending to infinity. 

We remark that as in the semi elliptic case (see \cite{NSV08} and especially the numerical experiments therein), if the equal functions $d_{\pm}$ have a root then we  expect an unpredictable asymptotical behavior of the eigenvalues of  coefficient   matrix $\cal{M_\alpha}$.

\subsubsection{Case III: General case }
In the case where $d_+\neq d_-$ the term $(2G_n-\mathbb{I}_{n})(T_n-T_{n}^{\mathrm{T}})$ is nonzero and it affects the spectrum of the preconditioned matrix. Specifically,
\begin{align}
&\tau_n^{-1}(p_{\alpha}(\theta))(\nu_{M,n}D_n^{-1}+T_n(p_{\alpha}(\theta))+(2G_n-\mathbb{I}_{n})T_n(2\mathbf{i}\Im\left\{g_{\alpha}(\theta)\right\}))\nonumber\\
&\hspace{1cm}\sim F_n^{-1/2}\mathbb{S}_n(\nu_{M,n}D_n^{-1}+T_n(p_{\alpha}(\theta))+(2G_n-\mathbb{I}_{n})T_n(2\mathbf{i}\Im\left\{g_{\alpha}(\theta)\right\}))\mathbb{S}_nF_n^{-1/2}\nonumber\\
&\hspace{1cm}=F_n^{-1/2}\mathbb{S}_n(\nu_{M,n}D_n^{-1})\mathbb{S}_nF_n^{-1/2}+F_n^{-1/2}\mathbb{S}_n(T_n(p_{\alpha}(\theta)))\mathbb{S}_nF_n^{-1/2}+ F_n^{-1/2}\mathbb{S}_n(2G_n-\mathbb{I}_{n})T_n(2\mathbf{i}\Im\left\{g_{\alpha}(\theta)\right\})\mathbb{S}_nF_n^{-1/2}, \nonumber
\end{align}
where  only the, new,  third term can add imaginary quantity on the eigenvalues.  However, we have observed through experimentation,   that the effect of this third  term on the real part of the eigenvalues is negligible. In this sense, we chose all  the  numerical experiments given in  Section 4  belong to this case mainly for  showing   the performance of our proposal  there were  our spectral analysis  do not explicitly and  in depth  cover the topic.

\subsection{Proposed Preconditioner: Two Dimensions}
In the two-dimensional case we use the second order spatial discretization, in order to be consistent with \cite{moghaderi171} and be able to readily compare the results. In this case, as reported in Section  \ref{sec:introduction:fde2d}, the coefficient matrix of the system is defined as
\begin{align}
\mathcal{M}_{(\alpha,\beta),N}^{(m)}&=\frac{1}{r}\mathbb{I}_{N}+D_+^{(m)}(\mathbb{I}_{n_2}\otimes S_{\alpha,n_1})+D_-^{(m)}(\mathbb{I}_{n_2}\otimes S_{\alpha,n_1}^{\mathrm{T}})+\frac{s}{r}\left(E_+^{(m)}( S_{\beta,n_2}\otimes \mathbb{I}_{n_1})+E_-^{(m)}( S_{\beta,n_2}^{\mathrm{T}}\otimes \mathbb{I}_{n_1})\right). \label{coef:2d}
\end{align}
 We recall that $S_{\alpha,n_1}=T_{n_1}(w_\alpha(\theta))$ and $S_{\beta,n_2}=T_{n_2}(w_\beta(\theta))$ (see \eqref{eq_w_function}, \eqref{S_as_Toeplitz}).
Again, for simplicity we here omit the time dependency in the notation.

Now let  $\mathcal{F}_{(\alpha,\beta)}(\theta_1,\theta_2)=q_\alpha(\theta_1)+\frac{s}{r}q_\beta(\theta_2)$ where $q$ is the real, nonnegative  and even function defined in \eqref{second_order_symbol}, $\theta_1, \theta_2 \in [-\pi,\pi],$ and $n_1$, $n_2$ the two integers   used for the discretization of  the domain $[L_x, R_x] \times [L_y, R_y]$.
Using the grid in \eqref{tau_preconditioner} we define the diagonal matrices
\begin{align}
F_{n_1,j}=&\mathrm{diag}(\mathcal{F}_{(\alpha,\beta)}(\theta_{i,n_1},\theta_{j,n_2}),  i=1,\ldots,n_1),\quad \nonumber
\end{align}
for each $j=1,\ldots,n_2$. Then, the $N\times N$ diagonal matrix is expressed as
\begin{align}
F_{N}=&
\begin{bmatrix*}
F_{n_1,1}&\\
&F_{n_1,2}&&\\
&&\ddots&\\
&&&\ddots&\\
&&&&F_{n_1,n_2}
\end{bmatrix*}. \label{F_N_matrix}
\end{align}
Let $\mathbb{S}_{n_1}$ and $\mathbb{S}_{n_2}$ be the discrete sine transform matrices of sizes $n_1$ and $n_2$, respectively,  as they defined in \eqref{sine_transform}.
Then, generalizing the idea of  (\ref{eq:prop_1D}),  our proposed preconditioner for this case is
\begin{align}\mathcal{P}_{\mathcal{F}_{(\alpha,\beta)},N}&=D_N\left(\mathbb{S}_{n_2} \otimes \mathbb{S}_{n_1}\right)F_N\left(\mathbb{S}_{n_2} \otimes \mathbb{S}_{n_1}\right), \label{2D_preconditioner}
\end{align}
where
\begin{align}
D_N&=(D_++D_-+E_++E_-)/4. \nonumber
\end{align}

The motivation of the above construction is to create a preconditioner that properly acts on the  different sources  affecting the spectrum of $\mathcal{M}_{(\alpha,\beta),N}.$ Specifically, the  diagonal part  operates   on the spatial space treating  the influence that  the coefficients of the equation have on the matrix, while the $\tau$ matrix  focuses  on the spectral space  and the ill-conditioning generated by the discretization  of the  fractional differential  operator. This observation  is a direct result of the GLT symbol associated to $\mathcal{M}_{(\alpha,\beta),N}$ and has been extensively studied  in \cite{NSV08} and \cite{V18}, for the case of semi elliptic differential equations. 
In the simplest, but not unusual  in applications, case where $d_{\pm}=d$, $e_{\pm}=e,$ we can counterbalance  the influence of the term $\frac{1}{r}$ in the spectrum of $\mathcal{M}_{(\alpha,\beta),N}$  incorporating  it  into  the $\tau$ part of the preconditioner.  Particularly,  we define  $\mathcal{\hat{F}}_{(\alpha,\beta)}(\theta_1,\theta_2)=\frac{1}{r}+d\cdot q_\alpha(\theta_1)+\frac{s}{r}e\cdot q_\beta(\theta_2)$ 
 replacing  the sampling of $\mathcal{F}_{(\alpha,\beta)}$ with that of $\mathcal{\hat{F}}_{(\alpha,\beta)}$ for the construction of $\hat{F}_N$ instead of $F_N$ in \eqref{F_N_matrix}.  Accordingly,  the new corresponding preconditioner $\mathcal{P}_{\mathcal{\hat{F}}_{(\alpha,\beta)},N}$  is defined as
 
\begin{align}
\mathcal{P}_{\mathcal{\hat{F}}_{(\alpha,\beta)},N}&=\left(\mathbb{S}_{n_2} \otimes \mathbb{S}_{n_1}\right)\hat{F}_N\left(\mathbb{S}_{n_2} \otimes \mathbb{S}_{n_1}\right). \label{2D_preconditioner_proof_case}
\end{align}

The following theorem shows that in this case,  the spectrum of the preconditioned matrix is bounded by positive constants independent of the size of the matrix.

\begin{theorem}\label{prop:main}
Assume that $d_{\pm}=d>0$, $e_{\pm}=e>0$. In this case the coefficient matrix of the system becomes

\begin{align}
A_N=\frac{1}{r}\mathbb{I}_N+(\mathbb{I}_{n_2} \otimes \hat{A}_{n_1}^\alpha)+(A_{n_2}^\beta \otimes \mathbb{I}_{n_1})=\left(\mathbb{I}_{n_2} \otimes (\frac{1}{r}\mathbb{I}_{n_1}+\hat{A}_{n_1}^\alpha)\right)+(A_{n_2}^\beta \otimes \mathbb{I}_{n_1})=\mathbb{I}_{n_2} \otimes A_{n_1}^\alpha+A_{n_2}^\beta \otimes \mathbb{I}_{n_1},
\end{align}
where
\begin{align}
A_{n_1}^\alpha&=\frac{1}{r}\mathbb{I}_N+T_{n_1}\left(d\cdot (w_\alpha(\theta)+w_\alpha(-\theta))\right)=T_{n_1}\left( \frac{1}{r}+d\cdot q_\alpha(\theta)\right),
\nonumber\\
A_{n_2}^\beta&= T_{n_2}\left(e\frac{s}{r}\cdot (w_\beta(\theta)+w_\beta(-\theta))\right)=T_{n_2}\left(e\frac{s}{r}\cdot q_\beta(\theta)\right). \nonumber
\end{align}
Then, the spectrum of the preconditioned matrix sequence $\left\{\mathcal{P}_{\mathcal{\hat{F}}_{(\alpha,\beta)},N}^{-1}A_N\right\}_N$ is bounded by  positive constants $c,C$ independent of  $N$.

\end{theorem}

\begin{proof}
We recall that
\begin{align}
h_x=(R_x-L_x)h_1,&\quad
h_y=(R_y-L_y)h_2,\nonumber\\
r=\frac{h_t}{2h_x^\alpha},&\quad s=\frac{h_t}{2h_y^\beta}, \nonumber
\end{align}
and 
\begin{align}
\hat{F}_N&=\mathbb{I}_{n_2}\otimes F_{n_1}^\alpha+F_{n_2}^\beta\otimes \mathbb{I}_{n_1}, \label{D}
\end{align}
where $\mathbb{I}_n$ is the identity matrix of order $n$ and
\begin{align}
F_{n_1}^\alpha&=\mathrm{diag}(d\cdot \mathcal{F}_\alpha(\theta_{i,n_1})+\frac{1}{r}),\quad i=1,\ldots,n_1, \label{h1} \\
F_{n_2}^\beta&=\mathrm{diag}(e\frac{s}{r}\cdot\mathcal{F}_\beta(\theta_{j,n_2})),\quad j=1,\ldots,n_2.\label{h2}
\end{align}
The matrix $A_N$ is  SPD since each of its terms is  a Kronecker product of a diagonal with a SPD Toeplitz matrix.
Hence,
\begin{align}
\mathcal{P}_{N}^{-1}A_{N} = \left(\mathbb{S}_{n_2} \otimes \mathbb{S}_{n_1}\right)\hat{F}_N^{-1}\left(\mathbb{S}_{n_2} \otimes \mathbb{S}_{n_1}\right)A_N,\nonumber
\end{align}
which is similar to the matrix
\begin{align}
 \hat{F}_N^{-1/2}\left(\mathbb{S}_{n_2} \otimes \mathbb{S}_{n_1}\right)A_N \left(\mathbb{S}_{n_2} \otimes \mathbb{S}_{n_1}\right)\hat{F}_N^{-1/2}.\nonumber
\end{align}
Thus,
\begin{align}
&\hat{F}_N^{-1/2}(S_{n_2} \otimes \mathbb{S}_{n_1})\left((\mathbb{I}_{n_2} \otimes A_{n_1}^\alpha) +(A_{n_2}^\beta \otimes \mathbb{I}_{n_1})\right)(\mathbb{S}_{n_2} \otimes \mathbb{S}_{n_1})\hat{F}_N^{-1/2}\nonumber\\
&\hspace{1cm}=\hat{F}_N^{-1/2}\left((\mathbb{S}_{n_2} \otimes \mathbb{S}_{n_1})(\mathbb{I}_{n_2} \otimes A_{n_1}^\alpha)(\mathbb{S}_{n_2} \otimes \mathbb{S}_{n_1})+(\mathbb{S}_{n_2} \otimes \mathbb{S}_{n_1})(A_{n2}^\beta \otimes \mathbb{I}_{n_1})(\mathbb{S}_{n_2} \otimes \mathbb{S}_{n_1})\right)\hat{F}_N^{-1/2}\nonumber\\
&\hspace{1cm}=\hat{F}_N^{-1/2}\left(\mathbb{I}_{n_2}\otimes \mathbb{S}_{n_1}A_{n_1}^\alpha S_{n_1} + \mathbb{S}_{n_2}A_{n_2}^\beta \mathbb{S}_{n_2} \otimes \mathbb{I}_{n_1}\right)\hat{F}_N^{-1/2}\nonumber\\
&\hspace{1cm}=\hat{F}_N^{-1/2}\left(\mathbb{I}_{n_2}\otimes (F_{n_1}^\alpha)^{1/2} (F_{n_1}^\alpha)^{-1/2} \mathbb{S}_{n_1}A_{n_1}^\alpha \mathbb{S}_{n_1}(F_{n_1}^\alpha)^{-1/2} (F_{n_1}^\alpha)^{1/2}\right. +\nonumber\\
&\hspace{2.13cm} +\left.(F_{n_2}^\beta)^{1/2} (F_{n_2}^\beta)^{-1/2}\mathbb{S}_{n_2}A_{n_2}^\beta \mathbb{S}_{n_2}(F_{n_2}^\beta)^{-1/2} (F_{n_2}^\beta)^{1/2} \otimes \mathbb{I}_{n_1}\vphantom{\frac{1}{r}}\right)\hat{F}_N^{-1/2}\nonumber\\
&\hspace{1cm}=\hat{F}_N^{-1/2}\left((\mathbb{I}_{n_2}\otimes (F_{n_1}^\alpha)^{1/2})\underbrace{(\mathbb{I}_{n_2}\otimes  (F_{n_1}^\alpha)^{-1/2} \mathbb{S}_{n_1}A_{n_1}^\alpha \mathbb{S}_{n_1}(F_{n_1}^\alpha)^{-1/2}}_{=L})(\mathbb{I}_{n_2}\otimes (F_{n_1}^\alpha)^{1/2})\right. + \nonumber\\
&\hspace{2.13cm} +\left.((F_{n_2}^\beta)^{1/2} \otimes \mathbb{I}_{n_1})\underbrace{((F_{n_2}^\beta)^{-1/2}\mathbb{S}_{n_2}A_{n_2}^\beta \mathbb{S}_{n_2}(F_{n_2}^\beta)^{-1/2}) \otimes \mathbb{I}_{n_1})}_{=R}((F_{n_2}^\beta)^{1/2} \otimes \mathbb{I}_{n_1})\vphantom{\frac{1}{r}}\right)\hat{F}_N^{-1/2}\nonumber\\
&\hspace{1cm}=\underbrace{\hat{F}_N^{-1/2}(\mathbb{I}_{n_2}\otimes (F_{n_1}^\alpha)^{1/2}) L (\mathbb{I}_{n_2}\otimes (F_{n_1}^\alpha)^{1/2})\hat{F}_N^{-1/2}}_{=A_L}+\underbrace{\hat{F}_N^{-1/2}((F_{n_2}^\beta)^{1/2} \otimes \mathbb{I}_{n_1})R((F_{n_2}^\beta)^{1/2} \otimes \mathbb{I}_{n_1})\hat{F}_N^{-1/2}}_{=A_R}. \label{precMatrix}
\end{align}
Let
\begin{align}
P_{n_1}^{\alpha}&=\mathbb{S}_{n_1}F_{n_1}^\alpha \mathbb{S}_{n_1},\nonumber\\
P_{n_2}^{\beta} &=\mathbb{S}_{n_2}F_{n_2}^\beta \mathbb{S}_{n_2}. \nonumber
\end{align}
Then, (see \cite{noutsos161}),  there exist positive  constants $c$ and $C$ independent of $n_1,n_2$, such that 
\begin{displaymath}
c<\sigma\left(\left(P_{n_1}^{\alpha}\right)^{-1}A_{n_1}^\alpha\right)<C \Rightarrow c<\sigma\left((F_{n_1}^\alpha)^{-1/2} \mathbb{S}_{n_1}A_{n_1}^\alpha \mathbb{S}_{n_1}(F_{n_1}^\alpha)^{-1/2}\right)<C,
\end{displaymath}
 and 
 \begin{displaymath}
 c<\sigma\left(\left(P_{n_2}^{\beta}\right)^{-1}A_{n_2}^\beta\right)<C\Rightarrow c<(F_{n_2}^\beta)^{-1/2}\mathbb{S}_{n_2}A_{n_2}^\beta \mathbb{S}_{n_2}(F_{n_2}^\beta)^{-1/2}<C.
\end{displaymath} 
  Consequently, for every normalized vector $x \in \mathbb{R}^N$ we find that: 
\begin{align}
c<x^{\mathrm{T}} L x<C, \qquad  c< x^{\mathrm{T}} R x<C. \nonumber
\end{align}
Since the  matrices $A_L,$ $A_R$ that form \eqref{precMatrix}  are SPD,  we  recall  some properties concerning such kind of matrices.  Specifically, we use the inequality $A>B$ for $A,B$ SPD matrices if $A-B >0$  is positive definite. In addition  if $A$, $B$, $C$, $D$, and $E$ are SPD, then
\begin{align}
A>B &\Leftrightarrow EAE>EBE, \label{eq:spd:1}\\
A>B \text{ and } C>D &\Leftrightarrow A+C>B+D.\label{eq:spd:2}
\end{align}
Therefore, we infer
\begin{align}
&\begin{cases}
c\mathbb{I}_N<L<C\mathbb{I}_N, \\
c\mathbb{I}_N<R<C\mathbb{I}_N,
\end{cases}\nonumber
\end{align}
and, using \eqref{eq:spd:1} and \eqref{eq:spd:2}, we deduce
\begin{align}
&\begin{cases}
 c\hat{F}_N^{-1}(\mathbb{I}_{n_2}\otimes F_{n_1}^\alpha)<A_L<C\hat{F}_N^{-1}(\mathbb{I}_{n_2}\otimes F_{n_1}^\alpha),\\
 c\hat{F}_N^{-1}(F_{n_2}^\beta \otimes \mathbb{I}_{n_1})<A_R<C\hat{F}_N^{-1}(F_{n_2}^\beta \otimes \mathbb{I}_{n_1}).
 \end{cases}\label{threeinequalities}
 \end{align}
Using again \eqref{eq:spd:1} and \eqref{eq:spd:2}, taking into account the two inequalities of \eqref{threeinequalities}, and \eqref{D}, we have
 \begin{align}
 c\hat{F}_N^{-1}(\mathbb{I}_{n_2}\otimes F_{n_1}^\alpha)+c\hat{F}_N^{-1}(F_{n_2}^\beta \otimes \mathbb{I}_{n_1}) &= c \hat{F}_N^{-1}\hat{F}_N = c\mathbb{I}_N,\nonumber\\
C\hat{F}_N^{-1}(\mathbb{I}_{n_2}\otimes F_{n_1}^\alpha)+C\hat{F}_N^{-1}(F_{n_2}^\beta \otimes \mathbb{I}_{n_1}) &= C \hat{F}_N^{-1}\hat{F}_N=C\mathbb{I}_N. \nonumber
\end{align}
Consequently, we conclude that 
\begin{displaymath} c\mathbb{I}_N\leq F_N^{-1/2}(\mathbb{S}_{n_1} \otimes \mathbb{S}_{n_2})A_N(\mathbb{S}_{n_1} \otimes \mathbb{S}_{n_2})F_N^{-1/2}\leq C \mathbb{I}_N. \end{displaymath}
 Therefore, the spectrum of the preconditioned matrix, which is similar to the $F_N^{-1/2}(\mathbb{S}_{n_1} \otimes \mathbb{S}_{n_2})A_N(\mathbb{S}_{n_1} \otimes \mathbb{S}_{n_2})F_N^{-1/2}$, lies in $[c,C]$. Moreover, from \cite{noutsos161} we expect all the eigenvalues to be clustered  at 1, something that is  numerically confirmed in the next section. 
 \end{proof}

 \begin{corollary}
Let the functions $d_{+}(x,y,t),d_{-}(x,y,t),e_{+}(x,y,t),e_{-}(x,y,t)$ being  strictly positive functions on $\Omega$, with $d_{+}(x,y,t)=d_{-}(x,y,t)=e_{+}(x,y,t)=e_{-}(x,y,t).$  Then,  the preconditioned matrix  sequence $\left\{\mathcal{P}_{\mathcal{\hat{F}}_{(\alpha,\beta)},N}^{-1} \mathcal{M}_{(\alpha,\beta),N}^{(m)} \right\}_N$ is bounded by  positive constants $c,C$ independent of  $N$.
\end{corollary}

\begin{proof}

The  proof  can be easily obtained  from the results of Theorem \ref{prop:main}  and the observation that the coefficient matrix in (\ref{coef:2d}) can be bounded by  

\begin{displaymath}
A^{c}_{N}\leq \mathcal{M}_{(\alpha,\beta),N}^{(m)} \leq A^{C}_N,
\end{displaymath}
where 
 \begin{align}
A_{N}^c&=\frac{1}{r}\mathbb{I}_{N}+c(\mathbb{I}_{n_2}\otimes S_{\alpha,n_1})+c(\mathbb{I}_{n_2}\otimes S_{\alpha,n_1}^{\mathrm{T}})+\frac{s\cdot c}{r}\left( (S_{\beta,n_2}\otimes \mathbb{I}_{n_1})+( S_{\beta,n_2}^{\mathrm{T}}\otimes \mathbb{I}_{n_1})\right), \nonumber
\end{align}
 \begin{align}
A_{N}^{C} &=\frac{1}{r}\mathbb{I}_{N}+C(\mathbb{I}_{n_2}\otimes S_{\alpha,n_1})+C (\mathbb{I}_{n_2}\otimes S_{\alpha,n_1}^{\mathrm{T}})+\frac{s\cdot C}{r}\left ( (S_{\beta,n_2}\otimes \mathbb{I}_{n_1})+( S_{\beta,n_2}^{\mathrm{T}}\otimes \mathbb{I}_{n_1})\right), \nonumber
\end{align}
and 
 \begin{displaymath}
c =\min_{(x,y,t)\in \Omega}\{d_{+}(x,y,t),d_{-}(x,y,t),e_{+}(x,y,t),e_{-}(x,y,t)\}, 
\end{displaymath}
\begin{displaymath}
 C=\max_{(x,y,t)\in \Omega}\{d_{+}(x,y,t),d_{-}(x,y,t), e_{+}(x,y,t),e_{-}(x,y,t)\}.
  \end{displaymath} 
Then, using Rayleigh quotient we obtain
\begin{displaymath}
\mathcal{P}_{\mathcal{\hat{F}}_N}^{-1} A^{c}_{N}  \leq \mathcal{P}_{\mathcal{\hat{F}}_N}^{-1} \mathcal{M}_{(\alpha,\beta),N}^{(m)} \leq \mathcal{P}_{\mathcal{\hat{F}}_N}^{-1} A^{C}_N
\end{displaymath}

\begin{displaymath}
\lambda_1(\mathcal{P}_{\mathcal{\hat{F}}_N}^{-1} A^{c}_{N})  \leq \lambda_1(\mathcal{P}_{\mathcal{\hat{F}}_N}^{-1} \mathcal{M}_{(\alpha,\beta),N}^{(m)})\leq \lambda_N(\mathcal{P}_{\mathcal{\hat{F}}_N}^{-1} \mathcal{M}_{(\alpha,\beta),N}^{(m)}) \leq \lambda_N(\mathcal{P}_{\mathcal{\hat{F}}_N}^{-1} A^{C}_N),
\end{displaymath}
and the proof is completed.
\end{proof}

\noindent In the subsequent section we report several numerical experiments which numerically confirm  that a similar spectral behavior of the preconditioned matrix is expected  in the more general case where  the coefficients functions of the equation are all different to each other.

\section{Numerical Examples}
\label{sec:numerical}
In this section we present three numerical examples to show the efficiency of the proposed preconditioners, compared with preconditioners discussed in~\cite{donatelli161} (one dimension) and \cite{moghaderi171} (two dimensions). We have chosen to compare our work with these works   since they are the most recent and have shown their superiority against the other proposals in the literature.  

\begin{itemize}
\item Example 1 is a one-dimensional problem, taken from~\cite[Example 1]{donatelli161}.  We compare and discuss the preconditioners therein with the proposed $\mathcal{P}_{\mathcal{F}_\alpha,n}$, and a few variations based on the spectral symbol. The fractional derivatives are of order $\alpha\in \{1.2,1.5,1.8\}$.

\item Example 2 is a two-dimensional problem, taken from~\cite[Example 1]{moghaderi171}.  We compare and discuss the preconditioners therein with the proposed $\mathcal{P}_{\mathcal{F}_{(\alpha,\beta)},N}$. The fractional derivatives are $\alpha=1.8$ and $\beta=1.6$.

\item Example 3 is the same experiment as Example 2, but with the fractional derivatives $\alpha=1.8$ and $\beta=1.2$.
\end{itemize}
The numerical experiments presented in Tables~\ref{tbl:table1}--\ref{tbl:table4} were implemented in \textsc{Julia} v1.1.0, using GMRES from the package \textsc{IterativeSolvers.jl} (GMRES tolerance is set to $10^{-7}$) and the \textsc{FFTW.jl} package. Benchmarking is done with \textsc{BenchmarkTools.jl} with 100 samplings and minimum time is presented in milliseconds. Experiments were run, in serial, on a computer with dual Intel Xeon E5 2630 v4 2.20 GHz (10 cores each) CPUs, and with 128 GB of RAM.

 Figures~\ref{fig:1a}, \ref{fig:1b}, \ref{fig:1c}, \ref{fig:2} and \ref{fig:3} show the scaled spectra of the preconditioned coefficient matrix $\mathcal{P}^{-1}\mathcal{M}_{\alpha,n_1}$ (and $\mathcal{P}^{-1}\mathcal{M}_{(\alpha,\beta),N}$) for different preconditioners $\mathcal{P}$, fractional derivatives $\alpha$, and matrix orders $n_1$ (and $\beta$, $N=n_1,n_2$). The scaling by a constant $c_0$ is performed as follows: We find the smallest disk enclosing all the eigenvalues of the considered matrix $A$. The center is denoted $c_0$ and the radius is $r$. Then, the spectrum is scaled as $\lambda_j(A)/c_0$ and the circle scaled and centered in $(1,0)$. The Julia package \textsc{BoundingSphere.jl} was used to compute $c_0$ and $r$ for all figures. The current scaling of the eigenvalues of preconditioned coefficient matrices is a visualization of the important effect  for the convergence rate of GMRES  of both the clustering and of the shape of the clustering.

In Tables~\ref{tbl:table1}--\ref{tbl:table4}, for each preconditioner, we present the number of iterations [it], minimal timing [ms], and the condition number of the preconditioned  matrix $\kappa$.  The best results are highlighted in bold.

\subsection{Example 1}
We compare the proposed preconditioner $\mathcal{P}_{\mathcal{F}_\alpha,n}$ with the ones presented in Example 1 from~\cite{donatelli161} (and two alternative symbol-based preconditioners).
We consider the one-dimensional form of~\eqref{eq:fde} in  the domain $[L_1,R_1]\times[t_0,T]= [0,2]\times [0,1]$,  where  the diffusion coefficients 
\begin{align}
d_+(x)&=\Gamma(3-\alpha)x^\alpha,\nonumber\\
d_-(x)&=\Gamma(3-\alpha)(2-x)^\alpha.\nonumber
\end{align}
are non-constant in space. Furthermore, the source term is
\begin{align}
f(x,t)=-32e^{-t}\left(x^2+\frac{(2-x)^2(8+x^2)}{8}-\frac{3(x^3+(2-x)^3)}{3-\alpha}+\frac{3(x^4+(2-x)^4)}{(4-\alpha)(3-\alpha)}\right),\nonumber
\end{align}
and the initial condition is
\begin{align}
u(x,0)=4x^2(2-x)^2,\nonumber
\end{align}
leading to the analytical solution $u(x,t)=4e^{-t}x^2(2-x)^2$. We assume $h_x=h_t=2/(n_1+1)$, that is, $\nu_{M,n_1}=h_x^{\alpha-1}$ and the number of time steps $M=(n_1+1)T/(R_1-L_1)=(n_1+1)/2$.
The set of fractional derivatives $\alpha$, for which a solution is computed for, is $\{1.2,1.5,1.8\}$ and in addition we consider the following set of partial dimensions for $n_1$: $\{2^6-1,2^7-1,2^8-1,2^9-1\}$.

In Table~\ref{tbl:table1} we present the results for the following preconditioners
\begin{itemize}
\item Identity ($\mathbb{I}_{n_1}$): GMRES without any preconditioner.
\item Circulant ($\mathcal{P}_{C,n_1}$): Described in~\cite{lei131} and implemented using FFT.
\item ``Full'' symbol ($\mathcal{P}_{\textsc{full},n_1}$): Defined as $\mathbb{S}_{n_1}\mathrm{diag}\left(\nu_{M,n_1}+d_{+,i}g_\alpha(\theta_{j,n_1})+d_{-,i}g_\alpha(-\theta_{j,n_1}), ~j=1,2,\ldots n_1\right)\mathbb{S}_{n_1}$ and implemented using FFT.
\item Symbol ($\mathcal{P}_{\mathcal{F}_\alpha,n_1}$): Proposed in Section~\ref{sec:prop1d}, $D_{n_1}\mathbb{S}_{n_1}\mathrm{diag}\left( p_\alpha(\theta_{j,n_1}), j=1,2,\ldots n_1\right)\mathbb{S}_{n_1}$, and implemented using FFT.
\end{itemize}

\begin{figure}[!ht]
\centering
\includegraphics[width=0.32\textwidth]{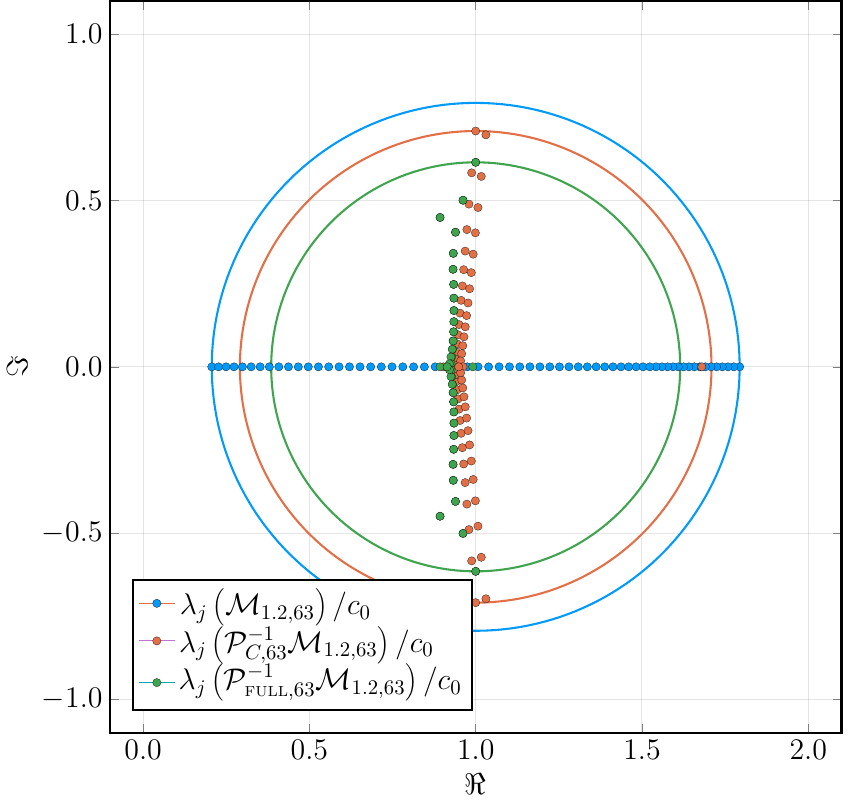}
\includegraphics[width=0.32\textwidth]{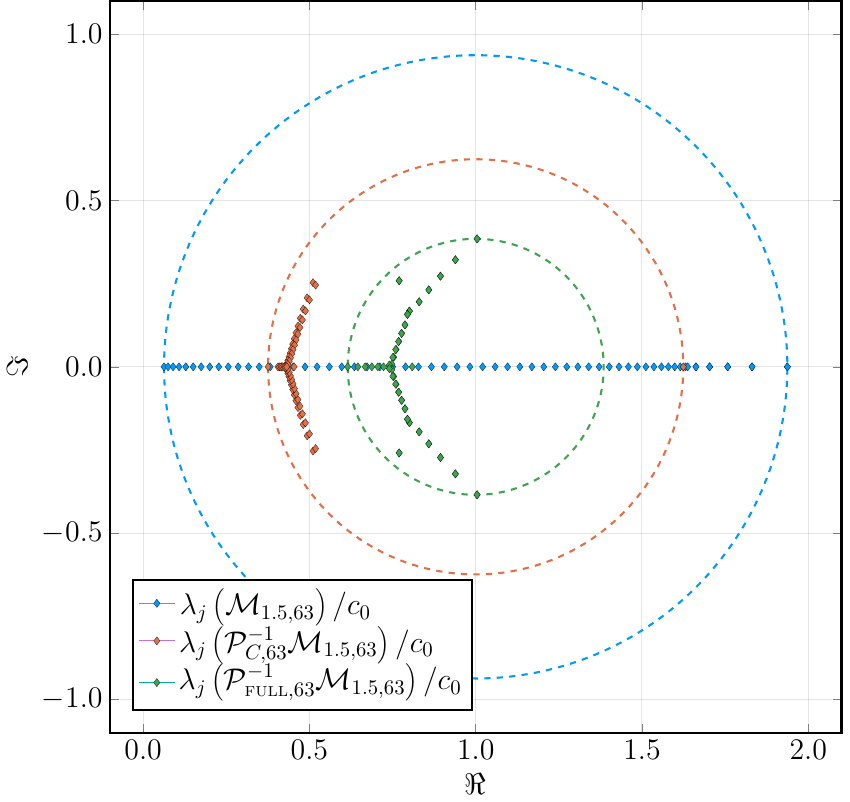}
\includegraphics[width=0.32\textwidth]{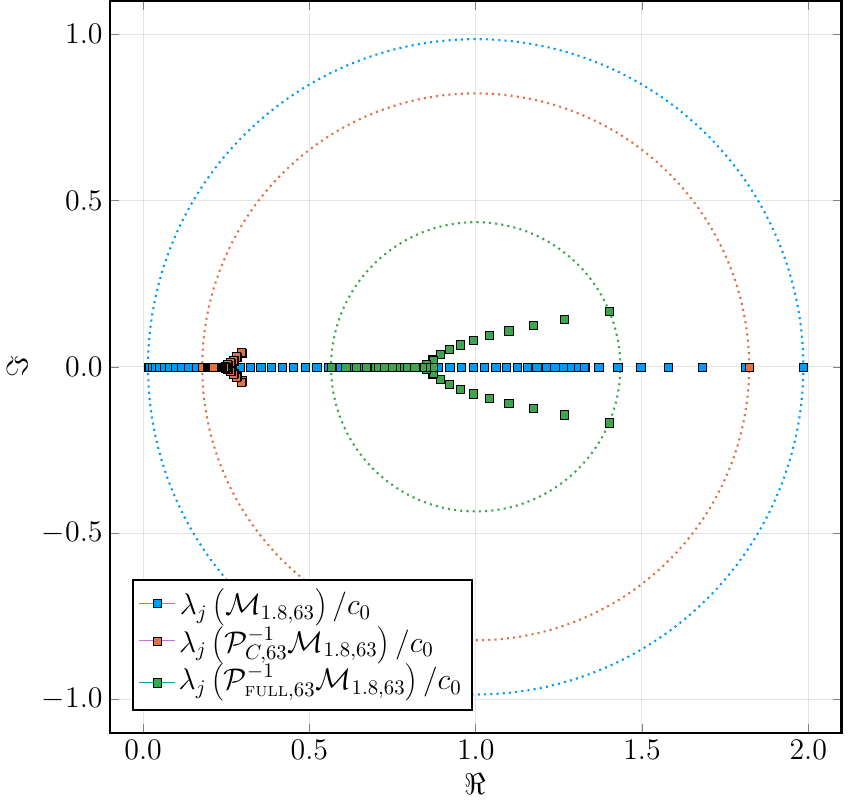}
\caption{[Example 1: 1D, $\alpha=\{1.2,1.5,1.8\}$] Scaled spectra of the resulting matrices when the preconditioners $\mathbb{I}_{n_1}$, $\mathcal{P}_{C,n_1}$, and
$\mathcal{P}_{\textsc{full},n_1}$
are applied to the coefficient matrices $\mathcal{M}_{\alpha,n_1}$ and $n_1=2^6-1$.
\textbf{Left:} $\alpha=1.2$.  \textbf{Middle:} $\alpha=1.5$. \textbf{Right:} $\alpha=1.8$.}
\label{fig:1a}
\end{figure}
In Figure~\ref{fig:1a} we present the scaled spectra of the resulting matrices, when the preconditioners $\mathbb{I}_{n_1}$, $\mathcal{P}_{C,n_1}$, and
$\mathcal{P}_{\textsc{full},n_1}$
are applied to the coefficient matrices $\mathcal{M}_{\alpha,n_1}$ when $n_1=2^6-1$ and $\alpha=1.2$ (left), $\alpha=1.5$ (middle), and $\alpha=1.8$ (right). We conclude that the spectral behavior resulting from the circulant and ``full'' symbol preconditioner resemble each other, but the condition number is lower for the ``full'' symbol preconditioner, as seen in Table~\ref{tbl:table1}.
\begin{figure}[!ht]
\centering
\includegraphics[width=0.48\textwidth]{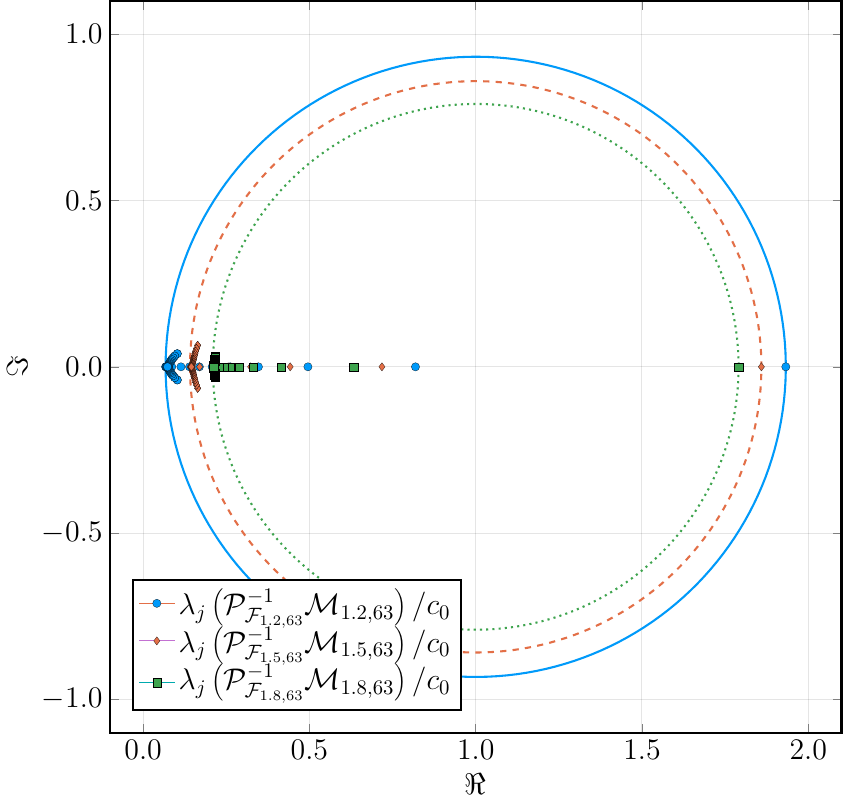}
\caption{[Example 1: 1D, $\alpha=\{1.2,1.5,1.8\}$] Scaled spectra of the resulting matrices when the preconditioners $\mathcal{P}_{\mathcal{F}_\alpha,n_1}$ are applied to the coefficient matrices $\mathcal{M}_{\alpha,n_1}$ for $n_1=2^6-1$.}
\label{fig:1b}
\end{figure}
In Figure~\ref{fig:1b} we show the scaled spectra of the resulting matrices when the preconditioners $\mathcal{P}_{\mathcal{F}_\alpha,n_1}$ are applied to the coefficient matrices $\mathcal{M}_{\alpha,n_1}$ with $n_1=2^6-1$ and $\alpha=\{1.2,1.5,1.8\}$.
We note that strong  clustering of the eigenvalues of the preconditioned matrices with  few large eigenvalues. The condition number is higher for the symbol preconditioner, compared to the ``full'' symbol preconditioner, however, as seen in Table~\ref{tbl:table1} both the number of iterations and execution time are lower for the symbol preconditioner. This numerically confirms what we  mentioned in Section 3.1 explaining the equation  (\ref{DM_analytical_presentation}), and the motivation of using a diagonal times a proper $\tau$ as preconditioner. In detail,  this  two terms preconditioner  properly acts on the different sources affecting the spectrum of the matrix: the diagonal part operates on the spatial space treating the influence that the coefficients of the equation have on the matrix, while the $\tau$ matrix focuses on the spectral space and the ill-conditioning generated by the discretization of the fractional differential operator.  Consequently,  this better clustering observed in Figure~\ref{fig:1b},  is the reason that the preconditioned GMRES method  (see \cite{axelsson861})   performs in general very well with  this  preconditioner.
In Table~\ref{tbl:table2} we present the results for the following preconditioners:
\begin{itemize}
\item First derivative ($\mathcal{P}_{1,n_1}$): Tridiagonal preconditioner based on the finite difference discretization of the first derivative, proposed in~\cite{donatelli161} and implemented using the Thomas algorithm.
\item Second derivative ($\mathcal{P}_{2,n_1}$): Tridiagonal preconditioner based on the finite difference discretization of the second derivative, proposed in~\cite{donatelli161} and implemented using the Thomas algorithm.
\item Tridiagonal ($\mathcal{P}_{\textsc{tri},n_1}$): Tridiagonal preconditioner based on the three main diagonals of the coefficient matrix and implemented using the Thomas algorithm.
\item Alternative symbol based ($\mathcal{P}_{\tilde{\mathcal{F}}_{\alpha},n_1}$): Constructed by $\mathbb{S}_{n_1}D_{n_1}\mathrm{diag}(p_\alpha(\theta_{j,n_1}))\mathbb{S}_{n_1}$ and implemented using FFT.
\end{itemize}
\begin{table}[!ht]
\centering
\caption{[Example 1: 1D, $\alpha=\{1.2,1.5,1.8\}$] Numerical experiments with GMRES and different preconditioners. For each preconditioner we present: average number of iterations for one time step [it], total timing in milliseconds [ms] to attain the approximate solution at time $T$, and the condition number $\kappa$ of the preconditioned  matrix, $\mathcal{P}^{-1}\mathcal{M}_{\alpha,n_1}$. The best results are highlighted in bold.}
\begin{tabular}{cc|ccc|ccc|ccc|ccccc}
\toprule
$\alpha$ & $n_1+1$&&$\mathbb{I}_{n_1}$&&&$\mathcal{P}_{\textsc{C},n_1}$&&&$\mathcal{P}_{\textsc{full},n_1}$&&&$\mathcal{P}_{\mathcal{F}_\alpha,n_1}$\\
&&\footnotesize [it]&\footnotesize [ms]&\footnotesize $\kappa$&\footnotesize [it]&\footnotesize [ms]&\footnotesize $\kappa$&\footnotesize [it]&\footnotesize [ms]&\footnotesize $\kappa$&\footnotesize [it]&\footnotesize [ms]&\footnotesize $\kappa$\\
\midrule
\footnotesize 1.2&\footnotesize $2^6$&\footnotesize \hfill 28.0&\footnotesize \hfill \textbf{1.7}&\footnotesize \hfill 9.6&\footnotesize \hfill 13.0&\footnotesize \hfill 9.6&\footnotesize \hfill 3.3&\footnotesize 14.0&\footnotesize \hfill 3.8&\footnotesize \textbf{1.6}&\footnotesize \hfill \textbf{7.2}&\footnotesize \hfill 2.3&\footnotesize \hfill 30.8\\
&\footnotesize $2^7$&\footnotesize \hfill 39.0&\footnotesize \hfill 24.3&\footnotesize \hfill 11.5&\footnotesize \hfill 14.0&\footnotesize \hfill 53.5&\footnotesize \hfill 3.6&\footnotesize 14.0&\footnotesize \hfill 17.6&\footnotesize \textbf{1.8}&\footnotesize \hfill \textbf{8.6}&\footnotesize \hfill \textbf{13.3}&\footnotesize \hfill 63.7\\
&\footnotesize $2^8$&\footnotesize \hfill 46.0&\footnotesize \hfill 114.9&\footnotesize \hfill 13.4&\footnotesize \hfill 13.0&\footnotesize \hfill 119.8&\footnotesize \hfill 3.8&\footnotesize 14.0&\footnotesize \hfill 68.8&\footnotesize \textbf{2.0}&\footnotesize \hfill \textbf{9.9}&\footnotesize \hfill \textbf{58.2}&\footnotesize \hfill 132.2\\
&\footnotesize $2^9$&\footnotesize \hfill 51.0&\footnotesize \hfill 594.5&\footnotesize \hfill 15.5&\footnotesize \hfill 12.0&\footnotesize \hfill 574.0&\footnotesize \hfill 4.2&\footnotesize 13.0&\footnotesize \hfill 312.7&\footnotesize \textbf{2.2}&\footnotesize \hfill \textbf{9.9}&\footnotesize \hfill \textbf{285.2}&\footnotesize \hfill 274.7\\
&\footnotesize $2^{10}$&\footnotesize \hfill 54.0&\footnotesize \hfill 2882.0 &\footnotesize \hfill 17.9&\footnotesize \hfill 11.0&\footnotesize \hfill 1927.0&\footnotesize \hfill 4.5&\footnotesize 12.0&\footnotesize \hfill \textbf{1415.0}&\footnotesize \hfill \textbf{2.4} &\footnotesize \hfill \textbf{10.9}&\footnotesize \hfill 1450.0&\footnotesize \hfill 571.4\\
&\footnotesize $2^{11}$&\footnotesize \hfill 56.0&\footnotesize \hfill  18569.0&\footnotesize \hfill 20.5&\footnotesize \hfill \textbf{10.0}&\footnotesize \hfill 11749.0&\footnotesize \hfill 4.9&\footnotesize 11.0&\footnotesize \hfill \textbf{8840.0}&\footnotesize \hfill \textbf{2.5}&\footnotesize \hfill 12.8&\footnotesize \hfill 9773.0&\footnotesize \hfill 1189.7 \\
\midrule
\footnotesize 1.5&\footnotesize $2^6$&\footnotesize \hfill 32.0&\footnotesize \hfill 2.0&\footnotesize \hfill 33.4&\footnotesize \hfill 12.0&\footnotesize \hfill 8.8&\footnotesize \hfill 7.1&\footnotesize 13.0&\footnotesize \hfill 3.2&\footnotesize \textbf{1.8}&\footnotesize \hfill \textbf{6.7}&\footnotesize \hfill \textbf{2.2}&\footnotesize \hfill 16.1\\
&\footnotesize $2^7$&\footnotesize \hfill 60.0&\footnotesize \hfill 37.2 &\footnotesize \hfill 51.2&\footnotesize \hfill 12.0&\footnotesize \hfill 46.7&\footnotesize \hfill 9.2&\footnotesize 13.0&\footnotesize \hfill 16.4&\footnotesize \textbf{2.1}&\footnotesize \hfill \textbf{8.0}&\footnotesize \hfill \textbf{12.5}&\footnotesize \hfill 33.3\\
&\footnotesize $2^8$&\footnotesize \hfill 89.0&\footnotesize \hfill 213.1&\footnotesize \hfill 75.8&\footnotesize \hfill 12.0&\footnotesize \hfill 111.3&\footnotesize \hfill 12.0&\footnotesize 13.0&\footnotesize \hfill 64.5&\footnotesize \textbf{2.3}&\footnotesize \hfill \textbf{8.5}&\footnotesize \hfill \textbf{52.6}&\footnotesize \hfill 70.9\\
&\footnotesize $2^9$&\footnotesize \hfill 122.0&\footnotesize \hfill 1389.0&\footnotesize \hfill 109.9&\footnotesize \hfill 12.0&\footnotesize \hfill 544.2&\footnotesize \hfill 15.8&\footnotesize 12.0&\footnotesize \hfill 288.9&\footnotesize \textbf{2.6}&\footnotesize \hfill \textbf{10.0}&\footnotesize \hfill \textbf{280.2}&\footnotesize \hfill 152.7\\
&\footnotesize $2^{10}$&\footnotesize \hfill 158.0&\footnotesize \hfill 8007.0&\footnotesize \hfill 157.7 &\footnotesize \hfill 11.0&\footnotesize \hfill 1779.0&\footnotesize \hfill 21.2&\footnotesize \hfill 11.0&\footnotesize \hfill \textbf{1366.0}&\footnotesize \hfill \textbf{2.9}&\footnotesize \hfill \textbf{10.0}&\footnotesize \hfill 1386.0&\footnotesize \hfill 331.8\\
&\footnotesize $2^{11}$&\footnotesize \hfill 195.0&\footnotesize \hfill 56266.0&\footnotesize \hfill 224.7&\footnotesize \hfill \textbf{10.0}&\footnotesize \hfill 11538.0&\footnotesize \hfill 28.6&\footnotesize \hfill \textbf{10.0} &\footnotesize \hfill \textbf{8551.0}&\footnotesize \hfill \textbf{3.2}&\footnotesize \hfill 11.0&\footnotesize \hfill 9142.0&\footnotesize \hfill 724.3\\
\midrule
\footnotesize 1.8&\footnotesize $2^6$&\footnotesize \hfill 32.0&\footnotesize \hfill 2.1&\footnotesize \hfill 136.5&\footnotesize \hfill 9.0&\footnotesize \hfill 6.6&\footnotesize \hfill 23.0&\footnotesize 10.0&\footnotesize \hfill 2.6&\footnotesize \textbf{2.6}&\footnotesize \hfill \textbf{6.1}&\footnotesize \hfill \textbf{2.2}&\footnotesize \hfill 9.7\\
&\footnotesize $2^7$&\footnotesize \hfill 67.0&\footnotesize \hfill 42.2&\footnotesize \hfill 266.3&\footnotesize \hfill 9.0&\footnotesize \hfill 36.1&\footnotesize \hfill 37.8&\footnotesize 11.0&\footnotesize \hfill 14.5&\footnotesize \textbf{2.8}&\footnotesize \hfill \textbf{6.8}&\footnotesize \hfill \textbf{11.2}&\footnotesize \hfill 19.5\\
&\footnotesize $2^8$&\footnotesize \hfill 131.0&\footnotesize \hfill 332.3&\footnotesize \hfill 494.8&\footnotesize \hfill 9.0&\footnotesize \hfill 89.8&\footnotesize \hfill 63.0&\footnotesize 10.0&\footnotesize \hfill 53.6&\footnotesize \textbf{2.9}&\footnotesize \hfill \textbf{7.0}&\footnotesize \hfill \textbf{47.2}&\footnotesize \hfill 40.8\\
&\footnotesize $2^9$&\footnotesize \hfill 231.2&\footnotesize \hfill 3085.0&\footnotesize \hfill 893.8&\footnotesize \hfill 9.0&\footnotesize \hfill 446.8&\footnotesize \hfill 106.3&\footnotesize \hfill 9.0&\footnotesize \hfill \textbf{257.9}&\footnotesize \textbf{2.9}&\footnotesize \hfill \textbf{8.6}&\footnotesize \hfill 262.8&\footnotesize \hfill 86.9\\
&\footnotesize $2^{10}$&\footnotesize \hfill 341.0&\footnotesize \hfill 20620.0&\footnotesize \hfill 1589.3&\footnotesize \hfill \textbf{8.0}&\footnotesize \hfill 1503.0&\footnotesize \hfill 180.5&\footnotesize \hfill \textbf{8.0}&\footnotesize \hfill \textbf{1191.0}&\footnotesize \hfill \textbf{3.0}&\footnotesize \hfill 10.0&\footnotesize \hfill 1370.0&\footnotesize \hfill 187.5\\
&\footnotesize $2^{11}$&\footnotesize \hfill 470.0&\footnotesize \hfill 163700.0&\footnotesize \hfill 2800.9&\footnotesize \hfill 8.0&\footnotesize \hfill 10197.0&\footnotesize \hfill 308.3&\footnotesize \hfill \textbf{7.0} &\footnotesize \hfill \textbf{7759.0}&\footnotesize \hfill \textbf{3.0}&\footnotesize \hfill 11.0&\footnotesize \hfill 9125.0&\footnotesize \hfill 408.1\\
\bottomrule
\end{tabular}
\label{tbl:table1}
\end{table}

\begin{table}[!ht]
\centering
\caption{[Example 1: 1D, $\alpha=\{1.2,1.5,1.8\}$] Numerical experiments with GMRES and different preconditioners. For each preconditioner we present the  average number of iterations for one time step [it], the total timing in milliseconds [ms] to attain the approximate solution at time $T$, and the condition number $\kappa$ of the preconditioned mass matrix, $\mathcal{P}^{-1}\mathcal{M}_{\alpha,n_1}$. The best results are highlighted in bold.}
\begin{tabular}{cc|ccc|ccc|ccc|cccccc}
\toprule
$\alpha$ & $n_1+1$&&$\mathcal{P}_{1,n_1}$&&&$\mathcal{P}_{2,n_1}$&&&$\mathcal{P}_{\textsc{tri},n_1}$&&&$\mathcal{P}_{\tilde{\mathcal{F}}_{\alpha},n_1}$\\
&&\footnotesize [it]&\footnotesize [ms]&\footnotesize $\kappa$&\footnotesize [it]&\footnotesize [ms]&\footnotesize $\kappa$&\footnotesize [it]&\footnotesize [ms]&\footnotesize $\kappa$&\footnotesize [it]&\footnotesize [ms]&\footnotesize $\kappa$\\
\midrule
\footnotesize 1.2&\footnotesize $2^6$&\footnotesize \hfill 8.0&\footnotesize \hfill 1.1&\footnotesize \hfill \textbf{1.2}&\footnotesize \hfill 9.0&\footnotesize \hfill 1.0&\footnotesize \hfill 2.1&\footnotesize \hfill \textbf{5.0}&\footnotesize \hfill \textbf{0.7}&\footnotesize \hfill 1.3&\footnotesize \hfill 7.5&\footnotesize \hfill 2.1&\footnotesize \hfill 29.2\\
&\footnotesize $2^7$&\footnotesize \hfill 8.0&\footnotesize \hfill 7.5&\footnotesize \hfill \textbf{1.3}&\footnotesize \hfill 10.0&\footnotesize \hfill 8.6&\footnotesize \hfill 2.2&\footnotesize \hfill \textbf{5.0}&\footnotesize \hfill \textbf{5.9}&\footnotesize \hfill 1.4&\footnotesize \hfill 8.5&\footnotesize \hfill 12.2&\footnotesize \hfill 58.7\\
&\footnotesize $2^8$&\footnotesize \hfill 7.0&\footnotesize \hfill 32.0&\footnotesize \hfill \textbf{1.3}&\footnotesize \hfill 10.0&\footnotesize \hfill 37.4&\footnotesize \hfill 2.4&\footnotesize \hfill \textbf{5.0}&\footnotesize \hfill \textbf{32.0}&\footnotesize \hfill 1.5&\footnotesize \hfill 9.9&\footnotesize \hfill 52.0&\footnotesize \hfill 118.6\\
&\footnotesize $2^9$&\footnotesize \hfill 7.0&\footnotesize \hfill 180.9&\footnotesize \hfill \textbf{1.4}&\footnotesize \hfill 10.0&\footnotesize \hfill 191.2&\footnotesize \hfill 2.6&\footnotesize \hfill \textbf{5.0}&\footnotesize \hfill \textbf{171.0}&\footnotesize \hfill 1.5&\footnotesize \hfill 9.9&\footnotesize \hfill 254.3&\footnotesize \hfill 239.7\\
&\footnotesize $2^{10}$&\footnotesize \hfill 6.0&\footnotesize \hfill 959.7 &\footnotesize \hfill \textbf{1.4}&\footnotesize \hfill 9.0&\footnotesize \hfill 1066.0&\footnotesize \hfill 2.8&\footnotesize \hfill \textbf{5.0}&\footnotesize \hfill \textbf{928.7}&\footnotesize \hfill 1.6&\footnotesize \hfill 11.0&\footnotesize \hfill 1363.0&\footnotesize \hfill 484.0\\
&\footnotesize $2^{11}$&\footnotesize \hfill 6.0&\footnotesize \hfill 7026.0&\footnotesize \hfill \textbf{1.5}&\footnotesize \hfill 9.0&\footnotesize \hfill 7675.0&\footnotesize \hfill 3.0&\footnotesize \hfill \textbf{5.0}&\footnotesize \hfill \textbf{6914.0}&\footnotesize \hfill 1.7&\footnotesize \hfill 12.0&\footnotesize \hfill 10787.0&\footnotesize \hfill 976.3\\
\midrule
\footnotesize 1.5&\footnotesize $2^6$&\footnotesize \hfill 16.0&\footnotesize \hfill 1.5&\footnotesize \hfill 2.5&\footnotesize \hfill 8.0&\footnotesize \hfill 1.0&\footnotesize \hfill \textbf{2.1}&\footnotesize \hfill \textbf{7.0}&\footnotesize \hfill \textbf{1.0}&\footnotesize \hfill 2.4&\footnotesize \hfill 8.7&\footnotesize \hfill 2.7&\footnotesize \hfill 13.6\\
&\footnotesize $2^7$&\footnotesize \hfill 20.0&\footnotesize \hfill 14.4&\footnotesize \hfill 3.1&\footnotesize \hfill 9.0&\footnotesize \hfill 8.1&\footnotesize \hfill \textbf{2.3}&\footnotesize \hfill \textbf{8.0}&\footnotesize \hfill \textbf{7.5}&\footnotesize \hfill 3.0&\footnotesize \hfill 8.0&\footnotesize \hfill 12.1&\footnotesize \hfill 26.3\\
&\footnotesize $2^8$&\footnotesize \hfill 24.0&\footnotesize \hfill 67.9&\footnotesize \hfill 4.0&\footnotesize \hfill 9.0&\footnotesize \hfill \textbf{35.3}&\footnotesize \hfill \textbf{2.7}&\footnotesize \hfill 11.0&\footnotesize \hfill 40.2&\footnotesize \hfill 4.0&\footnotesize \hfill \textbf{8.4}&\footnotesize \hfill 47.7&\footnotesize \hfill 51.8\\
&\footnotesize $2^9$&\footnotesize \hfill 26.0&\footnotesize \hfill 366.7&\footnotesize \hfill 5.2&\footnotesize \hfill 10.0&\footnotesize \hfill \textbf{197.5}&\footnotesize \hfill \textbf{3.0}&\footnotesize \hfill 13.0&\footnotesize \hfill 227.3&\footnotesize \hfill 5.4&\footnotesize \hfill \textbf{9.9}&\footnotesize \hfill 248.1&\footnotesize \hfill 103.0\\
&\footnotesize $2^{10}$&\footnotesize \hfill 27.0&\footnotesize \hfill 1810.0&\footnotesize \hfill 6.9&\footnotesize \hfill \textbf{10.0}&\footnotesize \hfill \textbf{1105.0}&\footnotesize \hfill \textbf{3.5}&\footnotesize \hfill 15.0&\footnotesize \hfill 1331.0&\footnotesize \hfill 7.4&\footnotesize \hfill \textbf{10.0}&\footnotesize \hfill 1636.0&\footnotesize \hfill 205.9\\
&\footnotesize $2^{11}$&\footnotesize \hfill 25.4&\footnotesize \hfill 11212.0&\footnotesize \hfill 9.0&\footnotesize \hfill \textbf{11.0}&\footnotesize \hfill \textbf{8179.0}&\footnotesize \hfill \textbf{4.0}&\footnotesize \hfill 18.0&\footnotesize \hfill 9684.0&\footnotesize \hfill 10.4&\footnotesize \hfill \textbf{11.0}&\footnotesize \hfill 10563.0&\footnotesize \hfill 424.5\\
\midrule
\footnotesize 1.8&\footnotesize $2^6$&\footnotesize \hfill 25.0&\footnotesize \hfill 2.5&\footnotesize \hfill 8.4&\footnotesize \hfill \textbf{6.0}&\footnotesize \hfill \textbf{0.8}&\footnotesize \hfill \textbf{1.6}&\footnotesize \hfill 7.0&\footnotesize \hfill 1.0&\footnotesize \hfill 3.5&\footnotesize \hfill 8.0&\footnotesize \hfill 2.3&\footnotesize \hfill 9.0\\
&\footnotesize $2^7$&\footnotesize \hfill 40.0&\footnotesize \hfill 27.3&\footnotesize \hfill 14.3&\footnotesize \hfill \textbf{6.0}&\footnotesize \hfill \textbf{6.3}&\footnotesize \hfill \textbf{1.7}&\footnotesize \hfill 10.0&\footnotesize \hfill 8.7&\footnotesize \hfill 5.6&\footnotesize \hfill 7.8&\footnotesize \hfill 11.3&\footnotesize \hfill 17.0\\
&\footnotesize $2^8$&\footnotesize \hfill 61.0&\footnotesize \hfill 159.8&\footnotesize \hfill 25.3&\footnotesize \hfill 7.0&\footnotesize \hfill \textbf{31.0}&\footnotesize \hfill \textbf{1.8}&\footnotesize \hfill 15.0&\footnotesize \hfill 48.3&\footnotesize \hfill 9.4&\footnotesize \hfill \textbf{6.9}&\footnotesize \hfill 43.3&\footnotesize \hfill 33.1\\
&\footnotesize $2^9$&\footnotesize \hfill 88.0&\footnotesize \hfill 1083.0&\footnotesize \hfill 44.7&\footnotesize \hfill \textbf{7.0}&\footnotesize \hfill \textbf{170.1}&\footnotesize \hfill \textbf{2.0}&\footnotesize \hfill 22.0&\footnotesize \hfill 325.4&\footnotesize \hfill 16.6&\footnotesize \hfill \textbf{7.0}&\footnotesize \hfill 222.4&\footnotesize \hfill 65.4\\
&\footnotesize $2^{10}$&\footnotesize \hfill 120.0&\footnotesize \hfill 6277.0&\footnotesize \hfill 78.8&\footnotesize \hfill \textbf{7.0}&\footnotesize \hfill \textbf{999.3}&\footnotesize \hfill \textbf{2.3}&\footnotesize \hfill 31.0&\footnotesize \hfill 1983.0&\footnotesize \hfill 30.0&\footnotesize \hfill 8.9&\footnotesize \hfill 1569.0&\footnotesize \hfill 130.1\\
&\footnotesize $2^{11}$&\footnotesize \hfill 158.0&\footnotesize \hfill 46716.0&\footnotesize \hfill 138.2&\footnotesize \hfill \textbf{7.0}&\footnotesize \hfill \textbf{7309.0}&\footnotesize \hfill \textbf{2.6}&\footnotesize \hfill 44.7&\footnotesize \hfill 15756.0&\footnotesize \hfill 54.6&\footnotesize \hfill 10.0&\footnotesize \hfill 10249.0&\footnotesize \hfill 259.8\\
\bottomrule
\end{tabular}
\label{tbl:table2}
\end{table}

\begin{figure}[!ht]
\centering
\includegraphics[width=0.32\textwidth]{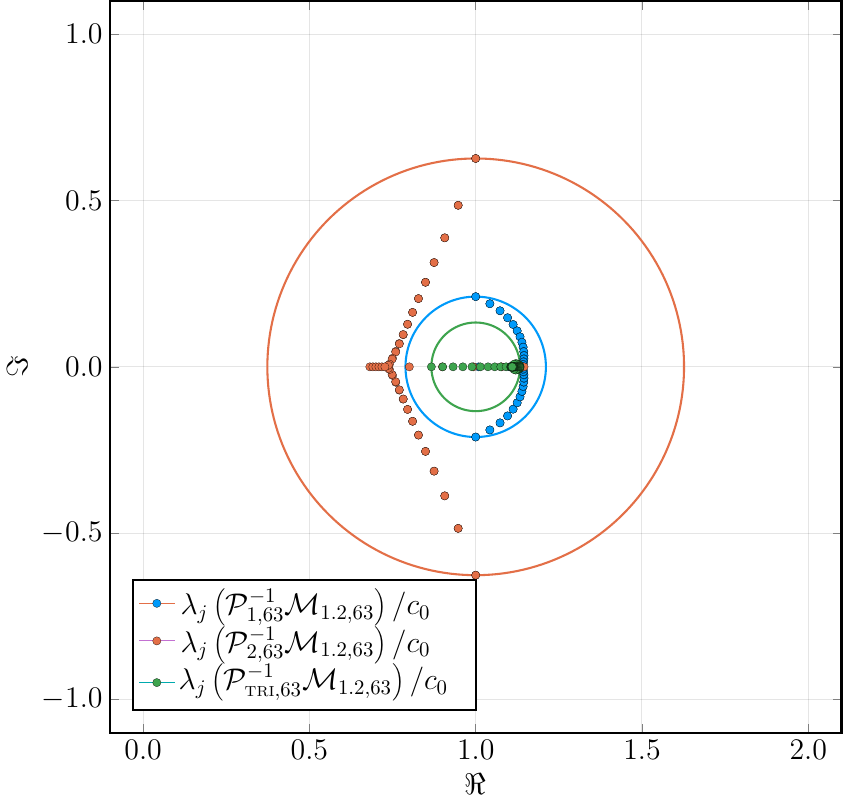}
\includegraphics[width=0.32\textwidth]{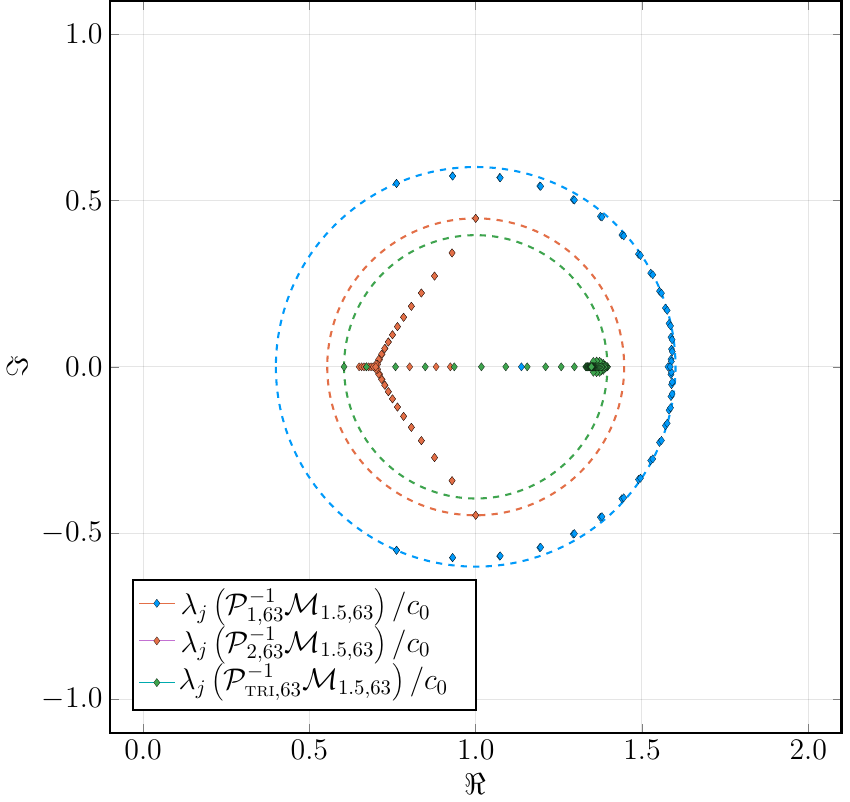}
\includegraphics[width=0.32\textwidth]{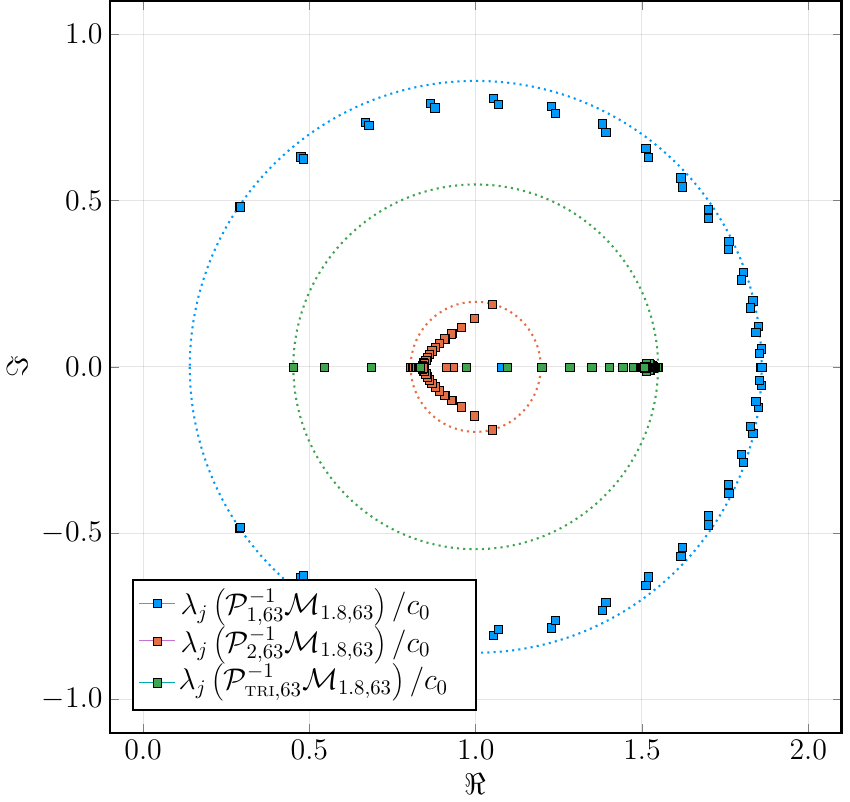}
\caption{[Example 1: 1D, $\alpha=\{1.2,1.5,1.8\}$] Scaled spectra of the resulting matrices when the preconditioners $\mathcal{P}_{1,n_1}$, $\mathcal{P}_{2,n_1}$, and $\mathcal{P}_{\textsc{tri},n_1}$ are applied to the  matrices $\mathcal{M}_{\alpha,n_1}$ and $n_1=2^6-1$.
\textbf{Left:} $\alpha=1.2$.  \textbf{Middle:} $\alpha=1.5$. \textbf{Right:} $\alpha=1.8$.}
\label{fig:1c}
\end{figure}

\begin{figure}[!ht]
\centering
\includegraphics[width=0.48\textwidth]{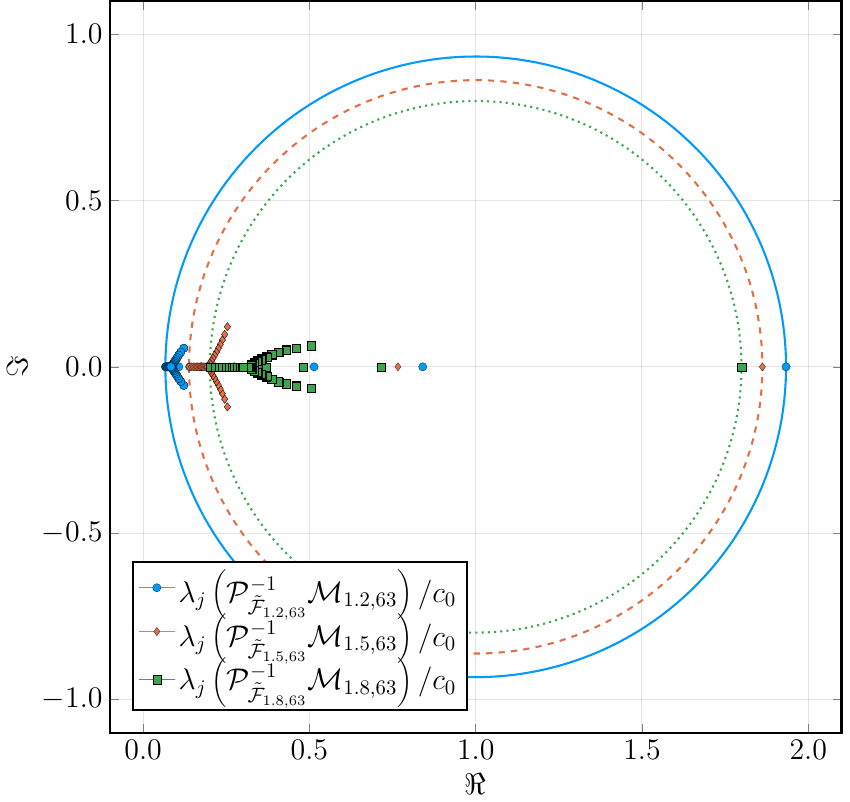}
\caption{[Example 1: 1D, $\alpha=\{1.2,1.5,1.8\}$] Scaled spectra of the resulting matrices when the preconditioners $\mathcal{P}_{\tilde{\mathcal{F}}_\alpha,n_1}$ are applied to the matrices $\mathcal{M}_{\alpha,n_1}$ for $n_1=2^6-1$.}
\label{fig:1d}
\end{figure}

As in Figure~\ref{fig:1a}, in Figure~\ref{fig:1c} we present the scaled spectra of the preconditioned matrix. The spectral behavior of the three preconditioners (first and second derivative and the tridiagonal) for different values of  $\alpha$ correlate well with the results presented in Table~\ref{tbl:table2}. In the left panel of Figure~\ref{fig:1c} the best clustering is obtained using the tridiagonal preconditioner, followed by the first derivative, and then by the second derivative. Since $\alpha=1.2$, a value close to one, this behavior is expected. When $\alpha=1.5$, as presented in the middle panel of Figure~\ref{fig:1c}, the results are similar for the three preconditioners, but the second derivative preconditioner performs in the best way  as $n_1$ increases. In the right panel of Figure~\ref{fig:1c} we see that the best clustering is observed for the second derivative preconditioner, and also show the best performances for all $n_1$ and all reported quantities (iterations, timings, and condition numbers).
The better performance of the preconditioners reported in Table~\ref{tbl:table2} as opposed the ones in Table~\ref{tbl:table1} is expected: this is due to the computational complexity of $\mathcal{O}(n)$ for the Thomas algorithm, as opposed to $\mathcal{O}(n\log n)$ for the DFT.        
In Figure~\ref{fig:1d} we present the scaled spectrum of an alternative symbol based preconditioner, $\mathcal{P}_{\tilde{\mathcal{F}}_\alpha,n_1}$, which performs slightly better than the proposed preconditioner $\mathcal{P}_{\mathcal{F}_\alpha,n_1}$ in Section~\ref{sec:prop1d} (compare Tables~\ref{tbl:table1} and \ref{tbl:table2}). This is mainly due to the avoided multiplication with the inverse of $D_n$ for $\mathcal{P}_{\tilde{\mathcal{F}}_\alpha,n_1}$, since the spectrum  of the resulted preconditioned matrices using  $\mathcal{P}_{\tilde{\mathcal{F}}_\alpha,n_1}$ and $\mathcal{P}_{\mathcal{F}_\alpha,n_1}$ are comparable. Furthermore, in this case it seems that the  most efficient choice of preconditioner is problem specific, depending on $d_{\pm}$.

\subsection{Example 2}
The considered two-dimensional example is originally from~\cite[Example 4.]{pang151} and is also discussed in~\cite[Example 1]{moghaderi171}. In \eqref{eq:fde}, define $\alpha=1.8$, $\beta=1.6$, and
\begin{align}
    d_+(x,y)=\Gamma(3-\alpha)(1+x)^\alpha(1+y)^2, \qquad d_-(x,y)=\Gamma(3-\alpha)(3-x)^\alpha(3-y)^2,\nonumber \\
    e_+(x,y)=\Gamma(3-\beta)(1+x)^2(1+y)^\beta,  \qquad e_-(x,y)=\Gamma(3-\beta)(3-x)^2(3-y)^\beta.\nonumber
\end{align}
The spatial domain is $\Omega=[0,2]\times[0,2]$ and the time interval is $[t_0,T]=[0,1]$.
The initial condition is
\begin{align}
    u(x,y,0)=u_0(x,y)=16x^2y^2(2-x)^2(2-y)^2,\nonumber
\end{align}
and the source term is
\begin{align}
f(x,y,t)&=-16e^{-t}\left(x^2(2-x)^2y^2(2-y)^2
+g_\alpha(x,y)
+g_\alpha(2-x,2-y)
+g_\beta(y,x)
+g_\beta(2-y,2-x)\right), \nonumber
\end{align} where 
\begin{align}
g_\gamma(x,y)&=\left(8x^{2-\gamma}-\frac{24x^{3-\gamma}}{3-\gamma}+\frac{24x^{4-\gamma}}{(4-\gamma)(3-\gamma)}\right)(1+x)^\gamma(1+y)^2y^2(2-y)^2,\nonumber
\end{align}
such that the solution to the FDE is given by $u(x,y,t)=16e^{-t}x^2(2-x)^2y^2(2-y)^2$.
Let $h=h_x=h_y=2/(n+1)$, with $n=n_1=n_2=M$, and $h_t=1/(M+1)$. Then,
\begin{align}
    \frac{1}{r}=\frac{2h^\alpha}{h_t}=\frac{2^{\alpha+1}M}{(n+1)^\alpha}=\frac{2^{\alpha+1}n}{(n+1)^\alpha}, \qquad \frac{s}{r}=\frac{h^\alpha}{h^\beta}=2^{\alpha-\beta}(n+1)^{\beta-\alpha}.\nonumber
\end{align}
In Table~\ref{tbl:table3} (and also Table~\ref{tbl:table4}) we present the results for the following preconditioners:
\begin{itemize}
\item Second derivative ($\mathcal{P}_{2,N}$): Preconditioner based on the finite difference discretization of the second derivative, proposed in~\cite{moghaderi171} and implemented using one Galerkin projection multigrid V-cycle.
\item Algebraic multigrid ($\mathcal{P}_{\textsc{MGM},N}$): Preconditioner based on algebraic multigrid, proposed in~\cite{moghaderi171} and implemented using one algebraic multigrid V-cycle.
\item Symbol ($\mathcal{P}_{\mathcal{F}_{(\alpha,\beta)},N}$): Proposed preconditioner and implemented using FFT.
\end{itemize}

We mention that in  multi-dimensional setting,   a negative results holds concerning the optimality of circulant algebra  when it is used for preconditioning  Toeplitz matrices  generated by  function with  zeros of order greater than one (e.g., see \cite{NSV_tcs,NSV_con}).  Thus, we find  a comparison with such kind of preconditioners to be  unnecessary . 
\begin{table}[!ht]
\centering
\caption{[Example 2: 2D, $\alpha=1.8, \beta=1.6$] Numerical experiments with GMRES and different preconditioners. For each preconditioner we present the average number of iterations for one time step [it], the total timing in milliseconds [ms] to attain the approximate solution at time $T$, and the condition number $\kappa$ of the preconditioned  matrix, $\mathcal{P}^{-1}\mathcal{M}_{(\alpha,\beta),N}$.  The best results are highlighted in bold.}
\begin{tabular}{c|ccc|ccc|ccc|cccc}
 \toprule
  $n_1=n_2$&&$\mathbb{I}_N$&&&$\mathcal{P}_{2,N}$&&&$\mathcal{P}_{\textsc{mgm},N}$&&&$\mathcal{P}_{\mathcal{F}_{(\alpha,\beta)},N}$ \\
  &\footnotesize [it]&\footnotesize [ms]&\footnotesize $\kappa$&\footnotesize [it]&\footnotesize [ms]&\footnotesize $\kappa$&\footnotesize [it]&\footnotesize [ms]&\footnotesize $\kappa$&\footnotesize [it]&\footnotesize [ms]&\footnotesize $\kappa$\\
 \midrule
 \footnotesize $2^4$&\footnotesize \hfill 37.0&\footnotesize \hfill 32.2&\footnotesize \hfill 57.4&\footnotesize \hfill 21.0&\footnotesize \hfill 64.8&\footnotesize \hfill 48.6&\footnotesize \hfill 10.0&\footnotesize \hfill 40.8&\footnotesize \hfill 3.7&\footnotesize \hfill \textbf{8.0}&\footnotesize \hfill \textbf{35.1}&\footnotesize \hfill \textbf{1.9}\\
 \footnotesize $2^5$&\footnotesize \hfill 73.0 &\footnotesize \hfill 331.4&\footnotesize \hfill 167.4&\footnotesize \hfill 17.6&\footnotesize \hfill 551.1&\footnotesize \hfill 31.7&\footnotesize \hfill 11.0&\footnotesize \hfill 383.1&\footnotesize \hfill 5.4&\footnotesize \textbf{8.0}&\footnotesize \hfill \textbf{296.8}&\footnotesize \hfill \textbf{2.7}\\
 \footnotesize $2^6$&\footnotesize \hfill 137.0 &\footnotesize \hfill 35440.0&\footnotesize \hfill 429.4&\footnotesize \hfill 17.0&\footnotesize \hfill 10465.0&\footnotesize \hfill 310.7&\footnotesize \hfill 11.0&\footnotesize \hfill 16146.0&\footnotesize \hfill 8.2&\footnotesize \hfill \textbf{9.0}&\footnotesize \hfill \textbf{6569.0}&\footnotesize \hfill \textbf{4.3}\\
 \footnotesize $2^7$&\footnotesize \hfill 251.0 &\footnotesize \hfill 1644134.0&\footnotesize \hfill 966.8&\footnotesize \hfill 17.0&\footnotesize \hfill 213713.0&\footnotesize \hfill 678.4&\footnotesize \hfill 10.0&\footnotesize 352471.0&\footnotesize \hfill 12.2&\footnotesize \hfill \textbf{9.0}&\footnotesize \hfill \textbf{135535.0}&\footnotesize \hfill \textbf{7.7}\\
 \bottomrule
\end{tabular}
\label{tbl:table3}
\end{table}

\begin{figure}[!ht]
\centering
\includegraphics[width=0.48\textwidth]{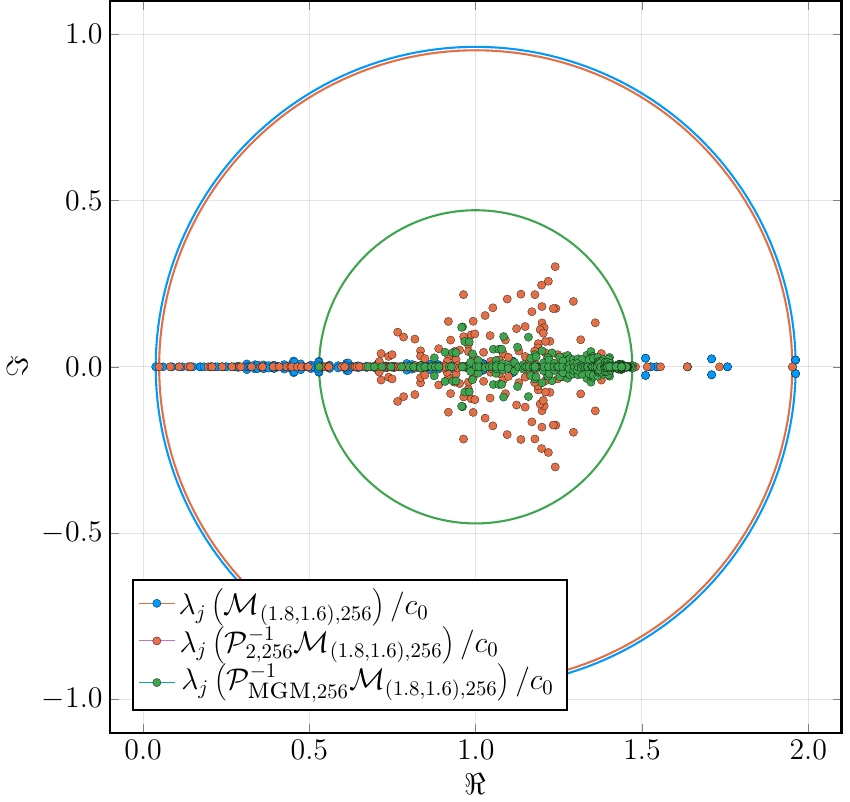}
\includegraphics[width=0.48\textwidth]{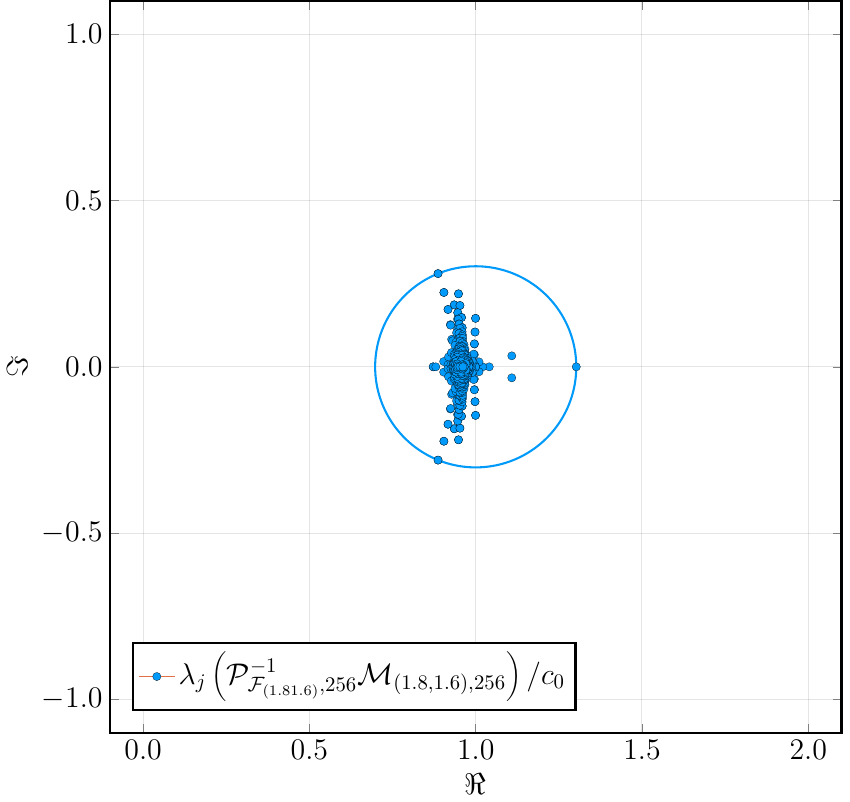}
\caption{[Example 2: 2D, $\alpha=1.8, \beta=1.6$] Scaled spectra of the resulting matrices when the preconditioners are applied to the  matrices $\mathcal{M}_{(\alpha,\beta),n_1^2}$ and $n_1=2^4$. \textbf{Left:} Preconditioners $\mathbb{I}_N$, $\mathcal{P}_{2,N}$, and $\mathcal{P}_{\textsc{mgm},N}$ \textbf{Right:} Preconditioner $\mathcal{P}_{\mathcal{F}_{(\alpha,\beta)},N}$.}
\label{fig:2}
\end{figure}

For details on the multigrid based preconditioners, $\mathcal{P}_{2,N}$ (Galerkin projection multigrid) and $\mathcal{P}_{\textsc{mgm},N}$ (algebraic multigrid), see~\cite{moghaderi171}.
The proposed symbol-based preconditioner, $\mathcal{P}_{\mathcal{F}_{(\alpha,\beta)},N}$, performs better than the multigrid-based preconditioners, as seen in Table~\ref{tbl:table3}. In Figure~\ref{fig:2} we present the scaled spectra of the preconditioned  matrices for $N=n_1n_2=2^8$. The clustering is better for the proposed symbol-based preconditioners than the other three, as seen comparing the left and right panels. We note in Table~\ref{tbl:table3} that the number of iterations are essentially constant both for the algebraic multigrid and the symbol-based preconditioners.

By fine tuning of the parameters for the multigrid-based preconditioners, such as number of smoothing steps, W-cycles etc, these results might be improved. However, the simplicity of the proposed preconditioner, where no fine-tunings are required, is advantageous.

\subsection{Example 3}
By modifying the coefficients $\alpha=1.8$ and $\beta=1.6$ in Example 2, to $\alpha=1.8$ and $\beta=1.2$ we obtain Example 3. In Table~\ref{tbl:table4} we present the same type of computations as in Table~\ref{tbl:table3}.
As discussed in~\cite{moghaderi171}, the performance of the proposed multigrid-based preconditioners depend on the fractional derivatives $\alpha$ and $\beta$. Since, in this example, $\alpha$ and $\beta$ differ more than in Example 2, and $\beta$ is far away from two, we clearly see in Table~\ref{tbl:table4} that the multigrid-based preconditioners perform worse than in Example 2. Especially note the worse behavior of the condition number for the algebraic multigrid-based preconditioner $\mathcal{P}_{\textsc{mgm},N}$. The condition numbers are essentially the same for the symbol-based preconditioner $\mathcal{P}_{\mathcal{F}_{(\alpha,\beta)},N}$ in Examples 2 and 3.

\begin{table}[!ht]
\centering
\caption{[Example 3: 2D, $\alpha=1.8, \beta=1.2$] Numerical experiments with GMRES and different preconditioners. For each preconditioner we present: average number of iterations for one time step [it], total timing in milliseconds [ms] to attain the approximate solution at time $T$, and the condition number $\kappa$ of the preconditioned  matrix, $\mathcal{P}^{-1}\mathcal{M}_{(\alpha,\beta),N}$. The best results are highlighted in bold.}
\begin{tabular}{c|ccc|ccc|ccc|cccc}
 \toprule
  $n_1=n_2$&&$\mathbb{I}_N$&&&$\mathcal{P}_{2,N}$&&&$\mathcal{P}_{\textsc{mgm},N}$&&&$\mathcal{P}_{\mathcal{F}_{(\alpha,\beta)},N}$ \\
  &\footnotesize [it]&\footnotesize [ms]&\footnotesize $\kappa$&\footnotesize [it]&\footnotesize [ms]&\footnotesize $\kappa$&\footnotesize [it]&\footnotesize [ms]&\footnotesize $\kappa$&\footnotesize [it]&\footnotesize [ms]&\footnotesize $\kappa$\\
 \midrule
 \footnotesize $2^4$ &\footnotesize \hfill 49.0 &\footnotesize \hfill 37.1 &\footnotesize \hfill 57.8&\footnotesize \hfill 26.5 &\footnotesize \hfill 79.7 &\footnotesize \hfill 42.8&\footnotesize \hfill 18.0 &\footnotesize \hfill 39.0 &\footnotesize \hfill 8.2&\footnotesize \hfill \textbf{10.0} &\footnotesize \hfill \textbf{37.0}&\footnotesize \hfill \textbf{1.9} \\
 \footnotesize $2^5$ &\footnotesize \hfill 92.0 &\footnotesize \hfill 394.0 &\footnotesize \hfill 162.9 &\footnotesize \hfill 32.0 &\footnotesize \hfill 713.8  &\footnotesize \hfill 104.0&\footnotesize \hfill 26.0 &\footnotesize \hfill 450.7 &\footnotesize \hfill 16.7&\footnotesize \textbf{12.0} &\footnotesize \hfill \textbf{329.0}&\footnotesize \hfill \textbf{2.7}    \\
 \footnotesize $2^6$ &\footnotesize \hfill  173.0 &\footnotesize \hfill 44532.0&\footnotesize \hfill  401.7 &\footnotesize \hfill 41.0 &\footnotesize \hfill 17197.0   &\footnotesize \hfill 231.6&\footnotesize \hfill 33.0 &\footnotesize \hfill 35021.0 &\footnotesize 32.8&\footnotesize \hfill \textbf{13.0} &\footnotesize \hfill \textbf{7493.0}&\footnotesize \hfill \textbf{4.4}   \\
 \footnotesize $2^7$ &\footnotesize  \hfill 316.0&\footnotesize \hfill 2070478.0&\footnotesize \hfill 876.4 &\footnotesize \hfill 51.0 &\footnotesize \hfill 438344.0 &\footnotesize \hfill 515.8&\footnotesize \hfill 41.0 &\footnotesize \hfill 1107711.0&\footnotesize 62.9&\footnotesize \hfill \textbf{14.5} &\footnotesize \hfill \textbf{171500.0}& \footnotesize \hfill \textbf{7.9}\\
 \bottomrule
\end{tabular}
\label{tbl:table4}
\end{table}

\begin{figure}[!ht]
\centering
\includegraphics[width=0.48\textwidth]{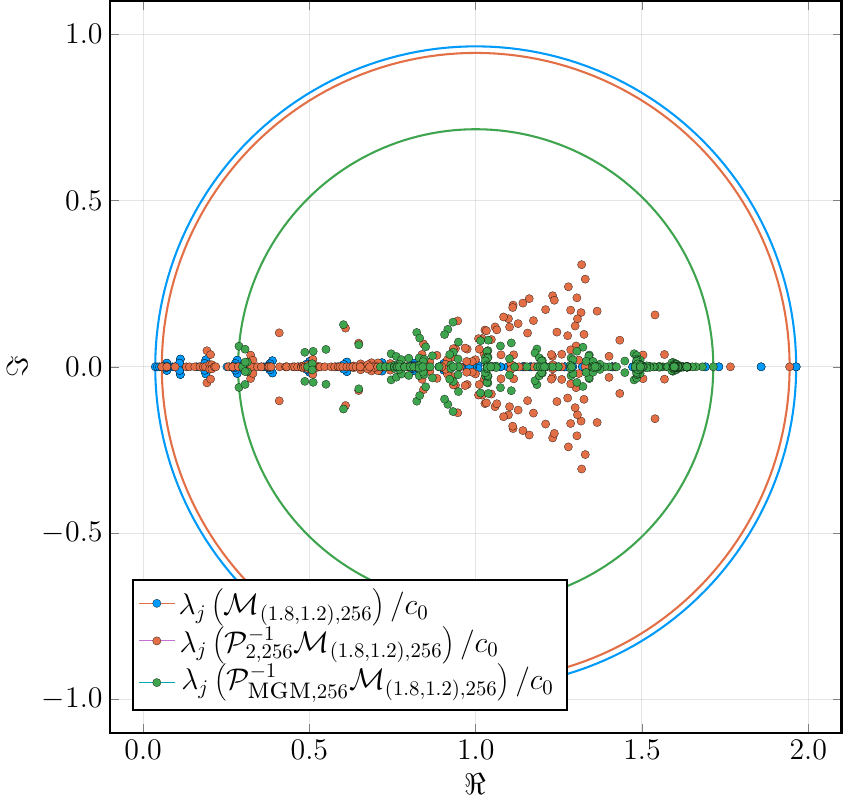}
\includegraphics[width=0.48\textwidth]{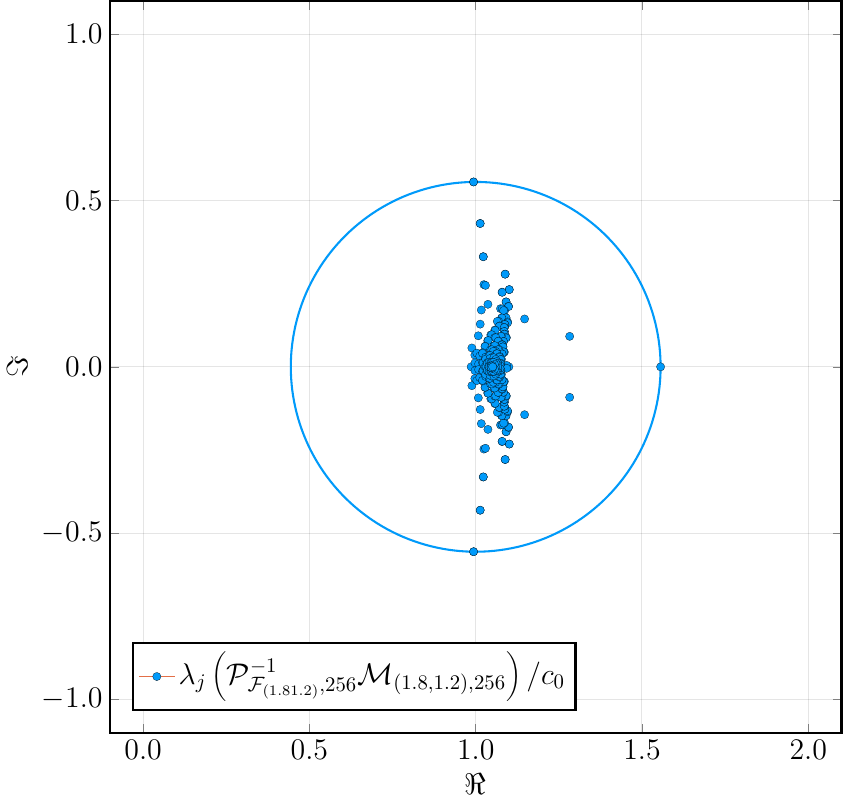}
\caption{[Example 3: 2D, $\alpha=1.8, \beta=1.2$] Scaled spectra of the resulting matrices when the preconditioners are applied to the coefficient  matrices $\mathcal{M}_{(\alpha,\beta),n_1^2}$, and $n_1=2^4$. \textbf{Left:} Preconditioners $\mathbb{I}_N$, $\mathcal{P}_{2,N}$, and $\mathcal{P}_{\textsc{mgm},N}$ \textbf{Right:} Preconditioner $\mathcal{P}_{\mathcal{F}_{(\alpha,\beta)},N}$.}
\label{fig:3}
\end{figure}

In Figure~\ref{fig:3} we present the same scaled spectra as in Figure~\ref{fig:2}, but regarding Example 3. Again, we note the advantageous clustering properties of the proposed symbol-based preconditioner in the right panel.

\section{Conclusions}
\label{sec:conclusions}
The purpose of the paper was the theoretical and numerical exploration of proper preconditioners based on the spectral symbols of the coefficient matrix for FDE problems. 
Beside the theoretical study, we have compared our results  with past ones, especially those presented in~\cite{donatelli161,moghaderi171}.
As expected, and  numerically shown in Example 1 which concerns  the one-dimensional case, our the proposed preconditioners performs slightly worse, at least in sequential computations, than the tridiagonal preconditions, because of the computational complexity involved. However, in the more challenging  two-dimensional case, as discussed in Examples 2 and 3, the proposed preconditioners do indeed perform better than the previously proposed multigrid-based preconditioners proposed and studied in \cite{moghaderi171}.

We note that future directions of research may include more complex problems, further analysis, and more extensive numerical experimentation. Also, problems where the fractional derivatives are close to three may be considered, since then we expect the symbol-based preconditioners to be even more advantageous, maybe even in the one-dimensional case.

\section{Acknowledgments}
The authors thank the anonymous reviewer's for suggestions improving the quality of the manuscript.
The second author was supported by the grant ``DRASI II'' of  Athens University of Economics and Business (Registration Number ER-3238).

\subsection*{CONFLICT OF INTEREST}
This study does not have any conflicts to disclose.

\bibliography{References2}
\bibliographystyle{wileyNJD-VANCOUVER}
\end{document}